\documentclass[sn-mathphys-num]{sn-jnl}% Math and Physical Sciences Numbered Reference Style 
\usepackage[dvipsnames,svgnames,x11names]{xcolor}
\usepackage{graphicx}%
\usepackage{multirow}%
\usepackage{amsmath,amssymb,amsfonts}%
\usepackage{amsthm}%
\usepackage{mathrsfs}%
\usepackage[title]{appendix}%
\usepackage{xcolor}%
\usepackage{textcomp}%
\usepackage{manyfoot}%
\usepackage{booktabs}%
\usepackage{algorithm}%
\usepackage{algorithmicx}%
\usepackage{algpseudocode}%
\usepackage{listings}%
\usepackage{geometry}
\usepackage{tikz}
\usetikzlibrary{patterns}
\usepackage{esint}
\tikzset{>=stealth}
\usepackage{bm}
\usepackage{yhmath}
\usepackage{makecell}
\numberwithin{equation}{section}
\theoremstyle{thmstyleone}%
\theoremstyle{thmstyletwo}%

\theoremstyle{thmstylethree}%

\begin{document}

\title[Article Title]{Complex variable solution on asymmetrical sequential shallow tunnelling in gravitational geomaterial considering static equilibrium}

%%=============================================================%%
%% GivenName	-> \fnm{Joergen W.}
%% Particle	-> \spfx{van der} -> surname prefix
%% FamilyName	-> \sur{Ploeg}
%% Suffix	-> \sfx{IV}
%% \author*[1,2]{\fnm{Joergen W.} \spfx{van der} \sur{Ploeg} 
%%  \sfx{IV}}\email{iauthor@gmail.com}
%%=============================================================%%

\author[1]{\fnm{Luo-bin} \sur{Lin}}\email{luobin\_lin@fjut.edu.cn}

\author*[2]{\fnm{Fu-quan} \sur{Chen}}\email{phdchen@fzu.edu.cn}
% \equalcont{These authors contributed equally to this work.}

\author[1]{\fnm{Chang-jie} \sur{Zheng}}\email{zcj@fjut.edu.cn}
% \equalcont{These authors contributed equally to this work.}

\author[1]{\fnm{Shang-shun} \sur{Lin}}\email{linshangshun@fjut.edu.cn}
% \equalcont{These authors contributed equally to this work.}

\affil[1]{\orgdiv{Fujian Provincial Key Laboratory of Advanced Technology and Informatization in Civil Engineering}, \orgname{Fujian University of Technology}, \orgaddress{\street{Xueyuan Road}, \city{Fuzhou}, \postcode{350118}, \state{Fujian}, \country{China}}}

\affil*[2]{\orgdiv{College of Civil Engineering}, \orgname{Fuzhou University}, \orgaddress{\street{Xueyuan Road}, \city{Fuzhou}, \postcode{350108}, \state{Fujian}, \country{China}}}

% \affil[3]{\orgdiv{Department}, \orgname{Organization}, \orgaddress{\street{Street}, \city{City}, \postcode{610101}, \state{State}, \country{Country}}}

%%==================================%%
%% Sample for unstructured abstract %%
%%==================================%%

\abstract{Asymmetrical sequential excavation is common in shallow tunnel engineering, especially for large-span tunnels. Owing to the lack of necessary conformal mappings, existing complex variable solutions on shallow tunnelling are only suitable for symmetrical cavities, and can not deal with asymmetrical sequential tunnelling effectively. This paper proposes a new complex variable solution on asymmetrical sequential shallow tunnelling by incorporating a bidirectional conformal mapping scheme consisting of Charge Simulation Method and Complex Dipole Simulation Method. Moreover, to eliminate the far-field displacement singularity of present complex variable method, a rigid static equilibrium mechanical model is established by fixing the far-field ground surface to equilibriate the nonzero resultant along cavity boundary due to graviational shallow tunnelling. The corresponding mixed boundary conditions along ground surface are transformed into homogenerous Riemann-Hilbert problems with extra constraints of traction along cavity boundaries, which are solved in an iterative manner to obtain reasonable stress and displacement fields of asymmetrical sequential shallow tunnelling. The proposed solution is validated by sufficient comparisons with equivalent finite element solution with good agreements. The comparisons also suggest that the proposed solution should be more accurate than the finite element one. A parametric investigation is finally conducted to illustrate possible practical applications of the proposed solution with several engineering recommendations. Additionally, the theoretical improvements and defects of the proposed solution are discussed for objectivity.}


%%\pacs[JEL Classification]{D8, H51}

%%\pacs[MSC Classification]{35A01, 65L10, 65L12, 65L20, 65L70}

\maketitle

\section{Introduction}
\label{sec:1}

\begin{figure}[htb]
  \centering
  \includegraphics[width = 0.6\textwidth]{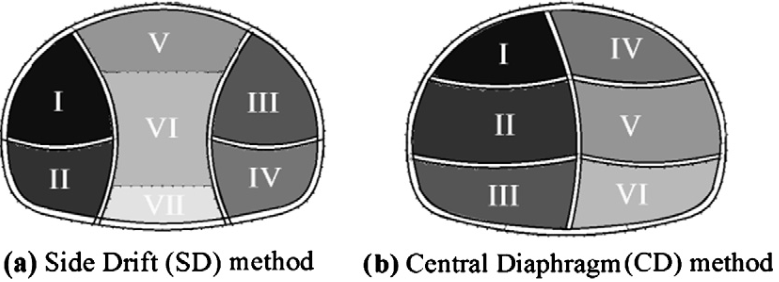}
  \caption{Typical sequential shallow excavation schemes with liners and temporary supports in Ref \cite{sharifzadeh13:_desig_niayes}}
  \label{fig:1}
\end{figure}

Sequential excavation is a widely used and effective construction method in large-span tunnels, for example, the Niayesh Road Tunnel of Niayesh and Sadr highways in Tehran \cite{sharifzadeh13:_desig_niayes}, the typical loess tunnel of the Zhengzhou-Xi'an high-speed railway \cite{li2016displacement}, the Hejie Tunnel of the Guiyang-Guangzhou high-speed railway \cite{fang2017shallow}, the Yingpan Road Tunnel in Changsha \cite{shi2017construction}, the Laodongshan Tunnel of the Guangtong-Kunming high-speed railway \cite{cao2018squeezing}, the Lianchengshan Tunnel of the Yinchuan-Kunming National Expressway \cite{liu2021deformation}, the Laohushan Tunnel of the Jinan Ring Expressway \cite{zhou2021stress}. According to excavation sequence, the construction methods can be classified into several catagories, such as side drift method, center diagram method, top heading and benche method, upper half vertical subdivision method, and so on. Fig. \ref{fig:1} shows typical side drift and center diagram schemes, and is cited from Ref \cite{sharifzadeh13:_desig_niayes}. 

These above sequential excavation methods are generally conducted by asymmetrical excavation transversely, and can be investigated using numerical methods (finite element method and finite difference method for instance). Though the numerical methods would deliver satisfactory results of sequential shallow tunnelling, several potential defects exist: (1) Long runtime of numerical models, especially for a full parametric analysis; (2) Multiple modellings for relavant sequential shallow tunnellings of different but similar geometries; (3) Insufficient insight of the mechanism of shallow tunnellings; and (4) Possible license deficiency of necessary commercial numerical softwares (ABAQUS and FLAC for instance). By contrast, analytical methods, such as complex variable method, can be used to adapt the above-mentioned potential defects of the numerical methods in the preliminary design stage of sequential shallow tunnelling.

The most classical complex variable solution on shallow tunnelling is the Verruijt solution \cite{Verruijt1997traction,Verruijt1997displacement}, in which the Verruijt conformal mapping is proposed to bidirectionally map a lower half plane containing a circular cavity onto a unit annulus. Based on the Verruijt conformal mapping, several extending solutions are developed to estimate the traction along circular cavity boundary \cite{zhang2018complex,wang2020analytical,zhang2020complex,zhang2021analytical}, to simulate the displacement along circular cavity boundary \cite{kong2019displacement,lu2019unified,kong2021analytical}, and to apply ground surface traction \cite{Wang2018surcharge_twin,Wang2018shallow_surcharge,Self2021APM_in_service,kong2021analytical}. However, no gravitational gradient of geomaterial is considered in the above extending solutions, which is a significant mechanical feature to distinguish shallow tunnelling from deep tunnelling. Strack \cite{Strack2002phdthesis} proposed the first complex variable solution on shallow tunnelling with consideration of gravitational gradient of geomaterial using Verruijt conformal mapping, and more relavant studies followed \cite{Strack_Verruijt2002buoyancy,Verruijt_Strack2008buoyancy,fang15:_compl,Lu2016,Lu2019new_solution}. The complex variable solutions mentioned above all focus on circular shallow tunnelling, while noncircular tunnels are also widely used in real-world engineering.

To adapt noncircular shallow tunnelling, Zeng et al. \cite{Zengguisen2019} proposed an extension of Verruijt conformal mapping by adding finite items in formation of bilateral Laurent series to backwardly map a unit annulus onto a lower half plane containing a noncircular but symmetrical cavity, and more complex variable solutions on noncircular shallow tunnelling by consideration of gravitational gradient are are proposed \cite{lu2021complex,zeng2022analytical,zeng2023analytical,fan2024analytical,zhou2024analytical}. The existing complex variable solutions on noncircular shallow tunnelling only focus on symmetrical and full cavity excavation, while asymmetrical and suquential excavation, which is commonly used in large-span tunnels, is rarely discussed. One possible reason is the lack of suitable conformal mapping. In this paper, we introduce a bidirectional conformal mapping scheme incorporating Charge Simulation Method \cite{amano1994charge} and Complex Dipole Simulation Method \cite{sakakibara2020bidirectional} to overcome such a mathematical obstable. Both simulation methods were originally put forward to solve electromagnetic problems, and were subsequently found efficient in constructing multiple types of conformal mappings.

Moreover, the above mentioned complex variable solutions considering gravitational gradient \cite{Zengguisen2019,lu2021complex,fan2024analytical,zhou2024analytical} generally ignore the violation of the very fundamental static equilibrium owing to the excavation of gravitational geomaterial. The nonzero resultant along cavity boundaries caused by unloading (buoyancy) of excavated geomaerial can not be equilibriated by any counter-acting force, and a displacement singularity would be consequently raised in far-field geomaterial. To establish a rigid static equilibrium, several recent solutions \cite{LIN2024appl_math_model,LIN2024comput_geotech_1,lin2024over-under-excavation,lin2024charge-simulation} have been put forward to construct new mechanical models by constraining far-field ground surface to artifically generate necessary counter-acting force to equilibriate the nonzero resultant caused by gravitational excavation. However, these recent solutions can not deal with asymmetrical sequential shallow tunnelling effectively.

In this paper, we seek a new complex variable solution to simultaneously solve both above mentioned problems: (i) Asymmetrical sequential shallow tunnelling with variety of cavity shapes, and (ii) Rigid static equilibrium to cancel far-field displacement singularity. With our new complex variable solution, reasonable stress and displacement fields of sequential shallow tunnelling with large-span cross section can be achieved, and the usage of complex variable method can further extended.

\section{Typical sequential shallow tunnelling}
\label{sec:2}

\subsection{Geometrical variations and notations}
\label{sec:2.1}

\begin{figure}[htb]
  \centering
  \includegraphics[width = 0.8\textwidth]{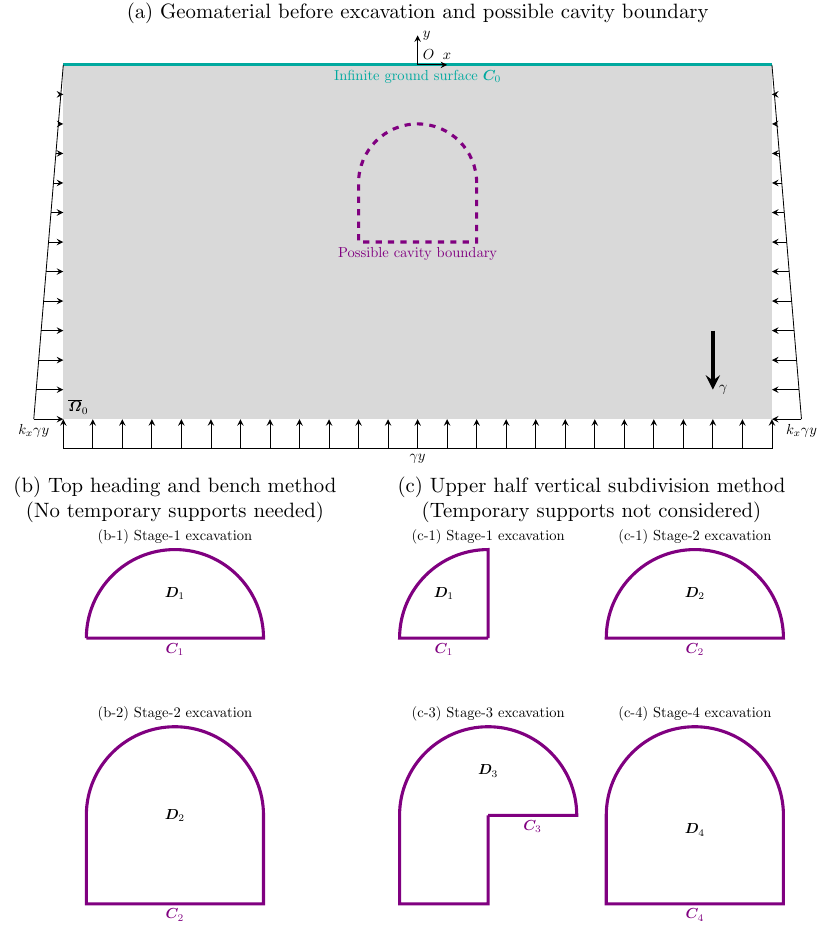}
  \caption{Geomaterial before excavation and two possible sequential exavation schemes}
  \label{fig:2}
\end{figure}

We may start from a gravitational geomaterial in Fig. \ref{fig:2}a. The ground surface is denoted by $ {\bm{C}}_{0} $, the geomaterial region is denoted by $ {\bm{\varOmega}}_{0} $, and the closure $ \overline{\bm{\varOmega}}_{0} = {\bm{C}}_{0} \cup {\bm{\varOmega}} $. The violet dash line denotes the possible final cavity boundary to be excavated. Two possible sequential shallow tunnelling schemes different from Fig. \ref{fig:1} are shown in Figs. \ref{fig:2}b and \ref{fig:2}c.

Fig. \ref{fig:2}b shows one possible excavation scheme using top heading and bench method, in which the cavity is excavated from top to bottom stepwisely, and no temporary supports may be needed. However, the symmetrically noncircular cavity geometry of top-heading-bench method in Fig. \ref{fig:2}b can already be mathematically studied by existing conformal mappings and complex variable solutions in above mentioned Refs \cite{Zengguisen2019,lu2021complex,zeng2022analytical,zeng2023analytical,fan2024analytical,zhou2024analytical}, and no obvious mathematical improvements of conformal mapping and complex variable method are thereby necessary.

By contrast, cavity geometry of the multi-stepwise upper half vertical subdivision method in Fig. \ref{fig:2}c would no longer hold symmetry, and the exsiting conformal mappings and complex variable solutions in Refs \cite{Zengguisen2019,lu2021complex,zeng2022analytical,zeng2023analytical,fan2024analytical,zhou2024analytical} would fail, and new conformal mapping scheme should be developed to form corresponding complex variable solution. The four-stage excavation scheme in Fig. \ref{fig:2}c is only a possible example of sequential shallow tunnelling for visualized illustration. Practical excavation schemes may be more complicated in real-world tunnel engineering, and temporary supports are generally needed. In this paper, we would focus on the mechanical variation of gravitational geomaterial alone caused by multi-stepwise upper half vertical subdivision method, and temporary supports are not considered.

To adapt the real-world complicated excavation schemes, we should use abstract notations to present possible excavation schemes in the following notations. A sequential shallow tunnelling is decomposed into $ N $ stages. The excavated region for $ j^{\rm{th}} $ stage ($ j = 1,2,3,\cdots,N $) are denoted by $ {\bm{D}}_{j} $, whose boundaries are denoted by by $ {\bm{C}}_{j} $. The above notations of regions and cavity boundaries are abstracted from Figs. \ref{fig:2}b and \ref{fig:2}c. After $ j^{\rm{th}} $ excavation, the remaining geomaterial region is reduced as a closure $ \overline{\bm{\varOmega}}_{j} = \overline{\bm{\varOmega}}_{0} - {\bm{D}}_{j} $. As shown in Fig. \ref{fig:2}c, for each stage, the excavated region after $ j^{\rm{th}} $ excavation always remains simply-connected, so that the rest geoamterial region $ \overline{\bm{\varOmega}}_{j} $ after $ j^{\rm{th}} $ excavation should always remain doubly-connected to hold the consistent connectivity. 

\subsection{Initial stress field and mechanical properties of geomaterial}
\label{sec:2.2}

As shown in Fig. \ref{fig:2}a, the geomaterial occupying the lower half plane $ \overline{\bm{\varOmega}}_{0} $ is assumed to be of small deformation with elasticity $ E $, Poisson's ratio $ \nu $, and shear modulus of $ G = \frac{E}{2(1+\nu)} $. A complex rectangular coordinate system $ z = x+{\rm{i}}y $ locating at some point of the ground surface $ {\bm{C}}_{0} $. An initial stress field of gravitational gradient $ \gamma $ and lateral coefficient $ k_{x} $ is orthotropically subjected through the geomaterial as
\begin{equation}
  \label{eq:2.1}
  \left\{
    \begin{aligned}
      & \sigma_{x}^{0}(z) = k_{x} \gamma y \\
      & \sigma_{y}^{0}(z) = \gamma y \\
      & \tau_{xy}^{0}(z) = 0 \\
    \end{aligned}
  \right.
  , \quad z = x + {\rm{i}}y \in \overline{\bm{\varOmega}}_{0}
\end{equation}
where $ \sigma_{x}^{0} $, $ \sigma_{y}^{0} $, and $ \tau_{xy}^{0} $ denote horizontal, vertical, and shear components of the initial stress field, respectively. Before excavation, the displacement in geomaterial is artifically reset to be zero.

\subsection{Sequential excavation and static equilibrium}
\label{sec:2.3}

Before any stage of sequential excavation, the traction along cavity boundary $ {\bm{C}}_{j} $ can be expressed as
\begin{equation}
  \label{eq:2.2}
  \left\{
    \begin{aligned}
      & X_{b}^{0,j}(S_{j}) = \sigma_{x}^{0}(S_{j}) \cdot \cos \langle \vec{n}_{j},\vec{x} \rangle + \tau_{xy}^{0}(S_{j}) \cdot \cos \langle \vec{n}_{j},\vec{y} \rangle \\
      & Y_{b}^{0,j}(S_{j}) = \tau_{xy}^{0}(S_{j}) \cdot \cos \langle \vec{n}_{j},\vec{x} \rangle + \sigma_{y}^{0}(S_{j}) \cdot \cos \langle \vec{n}_{j},\vec{y} \rangle \\
    \end{aligned}
  \right.
  , \quad S_{j} \in {\bm{C}}_{j}
\end{equation}
where $ X_{b}^{0,j} $ and $ Y_{b}^{0,j} $ denote horizontal and vertical components of surface traction of arbitrary point $ S_{j} $ along cavity boundary $ {\bm{C}}_{j} $ under the initial stress field from $ 0^{\rm{th}} $ stage to $ j^{\rm{th}} $ stage, respectively; $ \langle \vec{n}_{j},\vec{x} \rangle $ denotes the angle between outward normal $ \vec{n}_{j} $ and $ x $ axis, and $ \langle \vec{n}_{j},\vec{y} \rangle $ denotes the angle between outward normal $ \vec{n}_{j} $ and $ y $ axis.

The excavation after $ j^{\rm{th}} $ stage is conducted by mechanically applying opposite surface traction of Eq. (\ref{eq:2.2}) along cavity boundary $ {\bm{C}}_{j} $ in the integral formation as
\begin{equation}
  \label{eq:2.3}
  \int_{\bm{C}_{j}} \left[ X_{a}^{0,j}(S_{j}) + {\rm{i}} Y_{a}^{0,j}(S_{j}) \right]{\rm{d}}S_{j} = - \int_{\bm{C}_{j}} \left[ X_{b}^{0,j}(S_{j}) + {\rm{i}} Y_{b}^{0,j}(S_{j}) \right]{\rm{d}}S_{j}, \quad S_{j} \in {\bm{C}}_{j}
\end{equation}
where $ X_{a}^{0,j} $ and $ Y_{a}^{0,j} $ denote horizontal and vertical components of opposite surface traction of arbitrary point $ S_{j} $ along cavity boundary $ {\bm{C}}_{j} $ under the initial stress field from $ 0^{\rm{th}} $ stage to $ j^{\rm{th}} $ stage, respectively.

With the opposite surface traction in Eq. (\ref{eq:2.3}), the initial geomaterial $ \overline{\bm{\varOmega}}_{0} $ is reduced to $ \overline{\bm{\varOmega}}_{j} = \overline{\bm{\varOmega}}_{0} - {\bm{D}}_{j} $. The nonzero resultant of surface traction in Eq. (\ref{eq:2.3}) can be given as
\begin{equation}
  \label{eq:2.4}
  {\rm{i}}R_{y}^{0,j} = \varointclockwise_{\bm{C}_{j}} \left[ X_{a}^{0,j}(S_{j}) + {\rm{i}}Y_{a}^{0,j}(S_{j}) \right]{\rm{d}}S_{j} = {\rm{i}}\gamma \iint_{{\bm{D}}_{j}} {\rm{d}}x{\rm{d}}y
\end{equation}
where the resultant is located as arbitrary point $ z^{c}_{j} $ within cavity boundary $ {\bm{C}}_{j} $. It should be emphasized that the resultant $ {\rm{i}}R_{y}^{0,j} $ is always upward for $ j = 1,2,3,\cdots,N $ for whatever cavity shape in the sequential excavation stages (see Fig. \ref{fig:2}c-3), since the cavity boundary $ {\bm{C}}_{j} $ is closed and surrounds the simply-connected region $ {\bm{D}}_{j} $, which is crossed through by the downward potential gravitational gradient $ \gamma $. 

The nonzero resultant along cavity boundary $ {\bm{C}}_{j} $ would cause the remaining geomaterial of fully free ground surface to be a non-static equilibrium system, and a displacement singularity at infinity would be raised. Detailed mathematical proof procedure can be found in our previous study \cite{LIN2024appl_math_model}, and no repitition should be conducted in this paper.

\begin{figure}[htb]
  \centering
  \includegraphics[width = 0.8\textwidth]{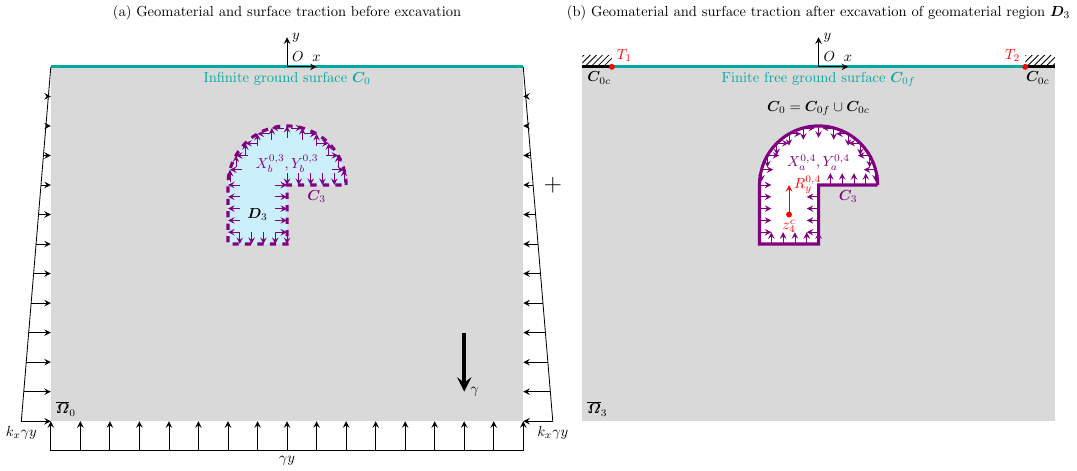}
  \caption{Boundary conditions of sequential excavation scheme of example in Fig. \ref{fig:2}c-3}
  \label{fig:3}
\end{figure}

To equilibriate the nonzero upward resultant along cavity boundary $ {\bm{C}}_{j} $, the far-field ground surface $ {\bm{C}}_{0c} $ is constrained to produce corresponding counter-acting force, and the remaining ground surface is left free and denoted by $ {\bm{C}}_{0f} $. In other words, $ {\bm{C}}_{0} = {\bm{C}}_{0c} \cup {\bm{C}}_{0f} $, as shown in Fig. \ref{fig:3}. As long as the width of segment $ {\bm{C}}_{0f} $ is large enough, the finite segment $ {\bm{C}}_{0f} $ can simulate an infinite ground surface. The following static equilibrium and mixed boundary conditions along ground surface can be established as
\begin{equation}
  \label{eq:2.5}
  \int_{\bm{C}_{0c}}\left[ X_{a}^{0,j}(T) + {\rm{i}}Y_{a}^{0,j}(T) \right]{\rm{d}}T = -{\rm{i}}R_{y}^{0,j}, \quad T \in {\bm{C}}_{0c}
\end{equation}
\begin{subequations}
  \label{eq:2.6}
  \begin{equation}
    \label{eq:2.6a}
    u_{c}^{0,j}(T) + {\rm{i}}v_{c}^{0,j}(T) = 0, \quad T \in {\bm{C}}_{0c}
  \end{equation}
  \begin{equation}
    \label{eq:2.6b}
    X_{f}^{0,j}(T) + {\rm{i}} Y_{f}^{0,j}(T) = 0, \quad T \in {\bm{C}}_{0f}
  \end{equation}
\end{subequations}
where $ u_{c}^{0,j} $ and $ v_{c}^{0,j} $ denote horizontal and vertical components of displacement along ground surface segment $ {\bm{C}}_{0c} $ due to the exacavation of geomaterial region $ {\bm{D}}_{j} $ from $ 0^{\rm{th}} $ stage to $ j^{\rm{th}} $ stage, respectively; $ X_{f}^{0,j} $ and $ Y_{f}^{0,j} $ dentote horizontal and vertical components of surface traction along ground surface segment $ {\bm{C}}_{0f} $ due to the excavation of geomaterial region $ {\bm{D}}_{j} $ from $ 0^{\rm{th}} $ stage to $ j^{\rm{th}} $ stage, respectively. Eqs. (\ref{eq:2.3}) and (\ref{eq:2.6}) form the necessary boundary conditions for excavation of $ {\bm{D}}_{j} $ geoamterial region from $ 0^{\rm{th}} $ stage to $ j^{\rm{th}} $ stage.

Fig. \ref{fig:3} present a corresponding case for graphic illustration of the abstract boundary conditions of sequential excavation in this section using the example of Stage-3 excavation in Fig. \ref{fig:2}c-3. In Fig. \ref{fig:3}a, the geomterial region $ {\bm{D}}_{3} $ is to be excavated, and the surface tractions along cavity boundary $ {\bm{C}}_{3} $ before excavation caused by the initial stress field is denoted by $ X_{b}^{0,3} $ and $ Y_{b}^{0,3} $. In Fig. \ref{fig:3}b, the geomaterial region $ {\bm{D}}_{3} $ is excavated, and the remaining geomaterial geometrically reduces from $ \overline{\bm{\varOmega}}_{0} $ to $ \overline{\bm{\varOmega}}_{3} = \overline{\bm{\varOmega}}_{0} - {\bm{D}}_{3} $. Meanwhile, the opposite surface tractions $ X_{a}^{0,3} $ and $ Y_{a}^{0,3} $ are applied along cavity boundary $ {\bm{C}}_{3} $. The upward resultant $ R_{y}^{0,3} $ due to excavation of geomaterial region $ {\bm{D}}_{3} $ is located at point $ z_{3}^{c} $ within cavity boundary $ {\bm{C}}_{3} $. The surface tractions $ X_{b}^{0,3} $ and $ Y_{b}^{0,3} $ in Fig. \ref{fig:3}a and the opposite surface traction $ X_{a}^{0,3} $ and $ Y_{a}^{0,3} $ in Fig. \ref{fig:3}b would cancel each other to free the boundary $ {\bm{C}}_{3} $.

With the mixed boundary conditions in Eqs. (\ref{eq:2.3}) and (\ref{eq:2.6}), the incremental stress and displacement fields in the remaining geomaterial can be expressed using the complex potentials as
\begin{subequations}
  \label{eq:2.7}
  \begin{equation}
    \label{eq:2.7a}
    \left\{
      \begin{aligned}
        & \sigma_{y}^{0,j}(z) + \sigma_{x}^{0,j}(z) = 2\left[ \varphi_{0,j}^{\prime}(z) + \overline{\varphi_{0,j}^{\prime}(z)} \right] \\
        & \sigma_{y}^{0,j}(z) - \sigma_{x}^{0,j}(z) + 2{\rm{i}}\tau_{xy}^{0,j}(z) = 2\left[ \overline{z}\varphi_{0,j}^{\prime\prime}(z) + \psi_{0,j}^{\prime}(z) \right] \\
      \end{aligned}
    \right.
    , \quad z \in \overline{\bm{\varOmega}}_{j}
  \end{equation}
  \begin{equation}
    \label{eq:2.7b}
    2G\left[ u_{0,j}(z) + {\rm{i}}v_{0,j}(z) \right] = \kappa \varphi_{0,j}(z) - z \overline{\varphi_{0,j}^{\prime}(z)} - \overline{\psi_{0,j}(z)}, \quad z \in \overline{\bm{\varOmega}}_{j}
  \end{equation}
\end{subequations}
where $ \sigma_{x}^{0,j} $, $ \sigma_{y}^{0,j} $, and $ \tau_{xy}^{0,j} $ denote horizontal, vertical, and shear components of incremental stress field due to full excavation of $ {\bm{D}}_{j} $ geomaterial region from $ 0^{\rm{th}} $ stage to $ j^{\rm{th}} $ stage, respectively; $ u_{0,j} $ and $ v_{0,j} $ denote horizontal and vertical components of incremental displacement field due to full excavation of $ {\bm{D}}_{j} $ geomaterial region from $ 0^{\rm{th}} $ stage to $ j^{\rm{th}} $ stage, respectively; $ \varphi_{0,j}^{\prime}(z) $ and $ \psi_{0,j}^{\prime}(z) $ denote the complex potentials to be determined using the mixed boundary conditions in Eqs. (\ref{eq:2.3}) and (\ref{eq:2.6}); $ \bullet^{\prime} $ and $ \bullet^{\prime\prime} $ denote the first- and second-order derivatives of abstract function $ \bullet $. For simplicity, the range of $ 0 \leq j \leq N $ is not repeated in the following deductions.

The total stress in geomaterial of $ j^{\rm{th}} $ excavation stage can be obtained by accumulation of the incremental one in Eq. (\ref{eq:2.7a}) and the initial one in Eq. (\ref{eq:2.1}) as
\begin{equation}
  \label{eq:2.7a'}
  \tag{2.7a'}
  \left\{
    \begin{aligned}
      & \sigma_{x}^{j}(z) = \sigma_{x}^{0}(z) + \sigma_{x}^{0,j}(z) \\
      & \sigma_{y}^{j}(z) = \sigma_{y}^{0}(z) + \sigma_{y}^{0,j}(z) \\
      & \tau_{xy}^{j}(z) = \tau_{xy}^{0}(z) + \tau_{xy}^{0,j}(z) \\
    \end{aligned}
  \right.
  , \quad z \in \overline{\bm{\varOmega}}_{j}, 0 \leq j \leq N
\end{equation}

\section{Riemann-Hilbert problem transformation}
\label{sec:3}

In this section, the complex variable method for sequential shallow tunnelling defined in Section \ref{sec:2.3} would be presented. To conduct the complex variable method, the physical geomaterial regions $ \overline{\bm{\varOmega}}_{j} $ ($ j = 1,2,3,\cdots,N $) should be bidirectionally mapped onto corresponding mapping unit annuli $ \overline{\bm{\omega}}_{j} $ via suitable conformal mappings, so that the boundary conditions, complex potentials, stres and displacement within the physical geomaterial regions $ \overline{\bm{\varOmega}}_{j} $ can be presented using the mathematically feasible Laurent series formations in the mapping unit annuli  $ \overline{\bm{\omega}}_{j} $ to facilitate further computation.

Our previous study \cite{lin2024charge-simulation} proposes a two-step conformal mapping scheme to bidirectionally map a lower half plane containing an arbitrary cavity onto a unit annulus. Each physical geomaterial region $ \overline{\bm{\varOmega}}_{j} $ accords with the geometrical and topological requirements of the the two-step conformal mapping in our previous study \cite{lin2024charge-simulation}, and can be thereby bidirectionally mapped onto corresponding mapping unit annuli $ \overline{\bm{\omega}}_{j} $ as
\begin{subequations}
  \label{eq:3.1}
  \begin{equation}
    \label{eq:3.1a}
    \zeta_{j} = \zeta_{j}(z) = \zeta_{j} \left[ w_{j} \left( z \right) \right], \quad z \in \overline{\bm{\varOmega}}_{j} \\
  \end{equation}
  \begin{equation}
    \label{eq:3.1b}
    z = z_{j}(\zeta_{j}) = z_{j} \left[ w_{j} \left( \zeta_{j} \right) \right], \quad \zeta_{j} \in \overline{\bm{\omega}}_{j} \\
  \end{equation}
\end{subequations}
The mapping scheme in Eq. (\ref{eq:3.1}) is briefly explained in Appendix \ref{sec:a}. Via the bidirectional conformal mapping, the whole geomaterial in the lower half plane $ \overline{\bm{\varOmega}}_{j} $ is bidirectionally mapped onto a bounded unit annulus $ \overline{\bm{\omega}}_{j} $; the cavity boundary $ {\bm{C}}_{j} $ is bidirectionally mapped onto the interior periphery of the unit annulus $ {\bm{c}}_{j} $; the constrained ground surface $ {\bm{C}}_{0c} $ and free one $ {\bm{C}}_{0f} $ are bidirectionally mapped onto corresponding segments $ {\bm{c}}_{0c,j} $ and $ {\bm{c}}_{0f,j} $ of exterior periphery of the unit annulus, respectively; the joint points $ T_{1} $ and $ T_{2} $ connecting $ {\bm{C}}_{0c} $ and $ {\bm{C}}_{0f} $ are bidirectionally mapped onto joint points $ t_{1,j} $ and $ t_{2,j} $ connecting $ {\bm{c}}_{0c,j} $ and $ {\bm{c}}_{0f,j} $, respectively. Joint points $ t_{1j} $ and $ t_{2,j} $ in mapping region $ \overline{\bm{\omega}}_{j} $ should remain singularities after conformal mapping. We should emphasize that $ \zeta_{j} = \rho_{j} \cdot \sigma_{j} = \rho_{j} \cdot {\rm{e}}^{{\rm{i}}\theta_{j}} $ holds in the following deductions.

With the conformal mapping scheme in Eq. (\ref{eq:3.1}), the complex potentials in Eq. (\ref{eq:2.7}) can be transformed as
\begin{equation}
  \label{eq:3.2}
  \left\{
    \begin{aligned}
      & \varphi_{0,j}(z) = \varphi_{0,j}(\zeta_{j}) \\
      & \psi_{0,j}(z) = \psi_{0,j}(\zeta_{j}) \\
    \end{aligned}
  \right.
  , \quad \zeta_{j} \in \overline{\bm{\omega}}_{j}
\end{equation}
Correspondingly, the incremental stress and displacement in Eq. (\ref{eq:2.7}) can be expressed using polar formation in mapping plane $ \zeta_{j} = \rho_{j} \cdot {\rm{e}}^{{\rm{i}}\theta_{j}} $ as
\begin{subequations}
  \label{eq:3.3}
  \begin{equation}
    \label{eq:3.3a}
    \sigma_{\theta}^{0,j}(\zeta_{j}) + \sigma_{\rho}^{0,j}(\zeta_{j}) = 2 \left[ \varPhi_{0,j}(\zeta_{j}) + \overline{\varPhi_{0,j}(\zeta_{j})} \right], \quad \zeta_{j} \in \overline{\bm{\omega}}_{j}
  \end{equation}
  \begin{equation}
    \label{eq:3.3b}
    \sigma_{\rho}^{0,j}(\zeta_{j}) + {\rm{i}} \tau_{\rho\theta}^{0,j}(\zeta_{j}) = \left[ \varPhi_{0,j}(\zeta_{j}) + \overline{\varPhi_{0,j}(\zeta_{j})} \right] - \frac{\overline{\zeta_{j}}}{\zeta_{j}} \left[ \frac{z_{j}(\zeta_{j})}{z_{j}^{\prime}(\zeta_{j})} \overline{\varPhi_{0,j}^{\prime}(\zeta_{j})} + \frac{\overline{z_{j}^{\prime}(\zeta_{j})}}{z_{j}^{\prime}(\zeta_{j})} \overline{\varPsi_{0,j}(\zeta_{j})} \right], \quad \zeta_{j} \in \overline{\bm{\omega}}_{j}
  \end{equation}
  \begin{equation}
    \label{eq:3.3c}
    g_{j}(\zeta_{j}) = 2G\left[ u_{0,j}(\zeta_{j}) + {\rm{i}} v_{0,j}(\zeta_{j}) \right] = \kappa \varphi_{0,j}(\zeta_{j}) - z_{j}(\zeta_{j})\overline{\varPhi_{0,j}(\zeta_{j})} - \overline{\psi_{0,j}(\zeta_{j})}, \quad \zeta_{j} \in \overline{\bm{\omega}}_{j}
  \end{equation}
  where $ \sigma_{\theta}^{0,j} $, $ \sigma_{\rho}^{0,j} $, and $ \tau_{\rho\theta}^{0,j} $ denote hoop, radial, and tangential components of the curvlinear stress mapped onto the mapping unit annuli $ \overline{\bm{\omega}}_{j} $, respectively; and
  \begin{equation*}
    \left\{
      \begin{aligned}
        & \varPhi_{0,j}(\zeta_{j}) = \frac{\varphi_{0,j}^{\prime}(\zeta_{j})}{z_{j}^{\prime}(\zeta_{j})} \\
        & \varPsi_{0,j}(\zeta_{j}) = \frac{\psi_{0,j}^{\prime}(\zeta_{j})}{z_{j}^{\prime}(\zeta_{j})} \\
      \end{aligned}
    \right.
    , \quad \zeta_{j} \in \overline{\bm{\omega}}_{j}
  \end{equation*}
\end{subequations}
The rectangular stress and displacement components in Eq. (\ref{eq:2.7}) can be computed using the curvilinear ones in Eq. (\ref{eq:3.3}) as
\begin{subequations}
  \label{eq:3.4}
  \begin{equation}
    \label{eq:3.4a}
    \left\{
      \begin{aligned}
        & \sigma_{y}^{0,j}(z) + \sigma_{x}^{0,j}(z) = \sigma_{\theta}^{0,j}\left[ \zeta_{j}(z) \right] + \sigma_{\rho}^{0,j} \left[ \zeta_{j}(z) \right] \\
        & \sigma_{y}^{0,j}(z) - \sigma_{x}^{0,j}(z) + 2{\rm{i}}\tau_{xy}^{0,j}(z) = \left\{ \sigma_{\theta}^{0,j}\left[ \zeta_{j}(z) \right] - \sigma_{\rho}^{0,j} \left[ \zeta_{j}(z) \right] + 2{\rm{i}} \tau_{\rho\theta}^{0,j} \left[ \zeta_{j}(z) \right] \right\} \cdot \frac{\overline{\zeta_{j}(z)}}{\zeta_{j}(z)} \cdot \frac{\overline{z_{j}^{\prime}[\zeta_{j}(z)]}}{z_{j}^{\prime}[\zeta_{j}(z)]}
      \end{aligned}
    \right.
    , \quad z \in \overline{\bm{\varOmega}}_{j}
  \end{equation}
  \begin{equation}
    \label{eq:3.4b}
    u_{0,j}(z) + {\rm{i}} v_{0,j}(z) = u_{0,j}\left[ \zeta_{j}(z) \right] + {\rm{i}} v_{0,j}\left[ \zeta_{j}(z) \right], \quad z \in \overline{\bm{\varOmega}}_{j}
  \end{equation}
\end{subequations}

With the conformal mapping scheme in Eq.(\ref{eq:3.1}), the mixed boundary conditions along ground surface in Eq. (\ref{eq:2.6}) can be mapped onto the mapping unit annuli $ \overline{\bm{\omega}}_{j} $ as
\begin{subequations}
  \label{eq:3.5}
  \begin{equation}
    \label{eq:3.5a}
    u_{c}^{0,j}(t_{j}) + {\rm{i}} v_{c}^{0,j}(t_{j}) = 0, \quad t_{j} \in {\bm{c}}_{0c,j}
  \end{equation}
  \begin{equation}
    \label{eq:3.5b}
    \sigma_{\rho}^{0,j}(t_{j}) + {\rm{i}} \tau_{\rho\theta}^{0,j}(t_{j}) = 0, \quad t_{j} \in {\bm{c}}_{0f,j}
  \end{equation}
\end{subequations}
where $ t_{j} = \sigma_{j} = {\rm{e}}^{{\rm{i}}\theta_{j}} $.

The segmental traction free boundary condition in Eq. (\ref{eq:3.5b}) can be used to establish analytic continuation and mixed boundary condition along ground surface regarding first-order derivative of complex potential $ \varPhi_{0,j}(\zeta_{j}) $. Substituting Eq. (\ref{eq:3.3b}) into Eq. (\ref{eq:3.5b}) yields
\begin{equation}
  \label{eq:3.6}
  \overline{\varPhi_{0,j}(t_{j})} = -\varPhi_{0,j}(t_{j}) + \frac{\overline{t}_{j}}{t_{j}}\frac{z_{j}(t_{j})}{z_{j}^{\prime}(t_{j})}\overline{\varPhi_{0,j}^{\prime}(t_{j})} + \frac{\overline{t}_{j}}{t_{j}}\frac{\overline{z_{j}^{\prime}(t_{j})}}{z_{j}^{\prime}(t_{j})}\overline{\varPsi_{0,j}(t_{j})}, \quad t_{j} \in {\bm{c}}_{0f,j}
\end{equation}
Partially substituting $ t_{j} = \overline{t}_{j}^{-1} $ into the items of the right-hand side of Eq. (\ref{eq:3.6}) yields
\begin{equation}
  \label{eq:3.7}
  \overline{\varPhi_{0,j}(t_{j})} = -\varPhi_{0,j}(\overline{t}_{j}^{-1}) + \frac{\overline{t}_{j}}{\overline{t}_{j}^{-1}}\frac{z_{j}(\overline{t}_{j}^{-1})}{z_{j}^{\prime}(\overline{t}_{j}^{-1})}\overline{\varPhi_{0,j}^{\prime}(t_{j})} + \frac{\overline{t}_{j}}{\overline{t}_{j}^{-1}}\frac{\overline{z_{j}^{\prime}(t_{j})}}{z_{j}^{\prime}(\overline{t}_{j}^{-1})}\overline{\varPsi_{0,j}(t_{j})}, \quad t_{j} \in {\bm{c}}_{0f,j}
\end{equation}
Replacing $ t_{j} = \sigma_{j} $ with $ \zeta_{j} = \rho_{j} \cdot \sigma_{j} $ in Eq. (\ref{eq:3.7}) yields
\begin{equation}
  \label{eq:3.8}
  \overline{\varPhi_{0,j}(\zeta_{j})} = -\varPhi_{0,j}(\overline{\zeta}_{j}^{-1}) + \overline{\zeta}_{j}^{2}\frac{z_{j}(\overline{\zeta}_{j}^{-1})}{z_{j}^{\prime}(\overline{\zeta}_{j}^{-1})}\overline{\varPhi_{0,j}^{\prime}(\zeta_{j})} + \overline{\zeta}_{j}^{2}\frac{\overline{z_{j}^{\prime}(\zeta_{j})}}{z_{j}^{\prime}(\overline{\zeta}_{j}^{-1})}\overline{\varPsi_{0,j}(\zeta_{j})}, \quad \zeta_{j} \in {\bm{\omega}}_{j}^{+}
\end{equation}
We should note that all items on the right-hand side of Eq. (\ref{eq:3.8}) are analytic in region $ \alpha_{j} \leq \rho_{j} < 1 $, which is contained within the exterior region $ \bm{\omega}_{j}^{-} $. Thus, Eq. (\ref{eq:3.8}) shows that $ \overline{\varPhi_{0,j}(\zeta_{j})} $ should be analytic in region $ \alpha_{j} \leq \rho_{j} < 1 $, indicating that analytic continuation has been conducted for $ \overline{\varPhi_{0,j}(\zeta_{j})} $, or simply its conjugate $ \varPhi_{0,j}(\zeta_{j}) $.

Substituting Eq. (\ref{eq:3.8}) into Eq. (\ref{eq:3.3b}) yields
\begin{equation}
  \label{eq:3.3b'}
  \tag{3.3b'}
  \begin{aligned}
    \sigma_{\rho}^{0,j}(\zeta_{j}) + {\rm{i}}\tau_{\rho\theta}^{0,j}(\zeta_{j}) = 
    & \; \overline{\zeta}_{j}^{2}\left[ \frac{z_{j}(\overline{\zeta}_{j}^{-1})}{z_{j}^{\prime}(\overline{\zeta}_{j}^{-1})} - \frac{1}{\zeta_{j}\overline{\zeta}_{j}}\frac{z_{j}(\zeta_{j})}{z_{j}^{\prime}(\zeta_{j})} \right]\overline{\varPhi_{0,j}^{\prime}(\zeta_{j})} + \overline{\zeta}_{j}^{2}\left[ \frac{\overline{z_{j}^{\prime}(\zeta_{j})}}{z_{j}^{\prime}(\overline{\zeta}_{j}^{-1})} - \frac{1}{\zeta_{j}\overline{\zeta}_{j}}\frac{\overline{z_{j}^{\prime}(\zeta_{j})}}{z_{j}^{\prime}(\zeta_{j})} \right]\overline{\varPsi_{j}(\zeta_{j})} \\
    & \; + \left[ \varPhi_{0,j}(\zeta_{j}) - \varPhi_{0,j}(\overline{\zeta}_{j}^{-1}) \right], \quad \zeta_{j} \in {\bm{\omega}}_{j}^{+}
  \end{aligned}
\end{equation}
The first-order of Eq. (\ref{eq:3.3c}) about $ \zeta_{j} $ can be expressed as
\begin{equation}
  \label{eq:3.9}
  \frac{{\rm{d}}g_{j}(\zeta_{j})}{{\rm{d}}\zeta_{j}} = \kappa z_{j}^{\prime}(\zeta_{j})\varPhi_{0,j}(\zeta_{j}) - z_{j}^{\prime}(\zeta_{j})\overline{\varPhi_{0,j}(\zeta_{j})} + z_{j}(\zeta_{j})\overline{\varPhi_{0,j}^{\prime}(\zeta_{j})} + \frac{\overline{\zeta}_{j}}{\zeta_{j}}\overline{z_{j}^{\prime}(\zeta_{j})}\;\overline{\varPsi_{0,j}(\zeta_{j})}
\end{equation}
Substituting Eq. (\ref{eq:3.8}) into Eq. (\ref{eq:3.9}) yields
\begin{equation}
  \label{eq:3.3c'}
  \tag{3.3c'}
  \begin{aligned}
    \frac{{\rm{d}}g_{j}(\zeta_{j})}{{\rm{d}}\zeta_{j}} = 
    & \; \frac{\overline{\zeta}_{j}}{\zeta_{j}}\left[ z_{j}(\zeta_{j}) - \zeta_{j}\overline{\zeta}_{j} \frac{z_{j}^{\prime}(\zeta_{j})}{z_{j}^{\prime}(\overline{\zeta}_{j}^{-1})}z(\overline{\zeta}_{j}^{-1}) \right]\overline{\varPhi_{0,j}^{\prime}(\zeta_{j})} + \frac{\overline{\zeta}_{j}}{\zeta_{j}}\left[ \overline{z_{j}^{\prime}(\zeta_{j})} - \zeta_{j}\overline{\zeta}_{j}\frac{z_{j}^{\prime}(\zeta_{j})}{z_{j}^{\prime}(\overline{\zeta}_{j}^{-1})}\overline{z_{j}^{\prime}(\zeta_{j})} \right]\overline{\varPsi_{0,j}(\zeta_{j})} \\
    & \; + \left[ \kappa z_{j}^{\prime}(\zeta_{j})\varPhi_{0,j}(\zeta_{j}) + z_{j}^{\prime}(\zeta_{j})\varPhi_{0,j}(\overline{\zeta}_{j}^{-1}) \right], \quad \zeta_{j} \in {\bm{\omega}}_{j}^{+}
  \end{aligned}
\end{equation}
We should notice that when $ \zeta_{j} \rightarrow \sigma_{j} $ (or $ \rho_{j} \rightarrow 1 $), the first lines of Eqs. (\ref{eq:3.3b'}) and (\ref{eq:3.3c'}) would turn to zero due to $ \overline{\sigma}_{j}^{-1} = \sigma_{j} $, and remaining second lines would form a homogenerous Riemann-Hilbert problem as \cite{LIN2024appl_math_model,lin2024over-under-excavation,lin2024charge-simulation}
\begin{subequations}
  \label{eq:3.10}
  \begin{equation}
    \label{eq:3.10a}
    \left. z_{j}^{\prime}(\zeta_{j})\left[ \sigma_{\rho}^{0,j}(\zeta_{j}) + {\rm{i}} \tau_{\rho\theta}^{0,j}(\zeta_{j}) \right] \right|_{\rho_{j} \rightarrow 1} = \varphi_{0,j}^{\prime+}(\sigma_{j}) - \varphi_{0,j}^{\prime -}(\sigma_{j}) = 0, \quad \sigma_{j} \in {\bm{c}}_{0f,j}
  \end{equation}
  \begin{equation}
    \label{eq:3.10b}
    \left. g_{j}^{\prime}(\zeta_{j}) \right|_{\rho_{j} \rightarrow 1} = \kappa \varphi_{0,j}^{\prime+}(\sigma_{j}) + \varphi_{0,j}^{\prime-}(\sigma_{j}) = 0, \quad \sigma_{j} \in {\bm{c}}_{0c,j}
  \end{equation}
\end{subequations}
where $ \varphi_{0,j}^{\prime+}(\sigma_{j}) $ and $ \varphi_{0,j}^{\prime-}(\sigma_{j}) $ denote the boundary values of $ \varphi_{0,j}^{\prime}(\zeta_{j})=z_{j}^{\prime}(\zeta_{j})\varPhi_{0,j}(\zeta_{j}) $ approaching boundary $ {\bm{c}}_{0,j} $ from $ {\bm{\omega}}^{+}_{j} $ and $ {\bm{\omega}}_{j}^{-} $ sides, respectively. 

The integrand of the left-hand side of Eq. (\ref{eq:2.3}) can be mapped from the physical plane $ z = x+{\rm{i}}y $ onto corresponding mapping plane $ \zeta_{j} = \rho_{j} \cdot {\rm{e}}^{{\rm{i}}\theta_{j}} $ according to the backward conformal mapping in Eq. (\ref{eq:3.1b}) as
\begin{subequations}
  \label{eq:3.11}
  \begin{equation}
    \label{eq:3.11a}
    X_{a}^{0,j}(S_{j}) + {\rm{i}} Y_{a}^{0,j}(S_{j}) = \frac{z_{j}^{\prime}(s_{j})}{|z_{j}^{\prime}(s_{j})|} \frac{{\rm{d}}s_{j}}{|{\rm{d}}s_{j}|} \cdot \left[ \sigma_{\rho}^{0,j}(s_{j}) + {\rm{i}} \tau_{\rho\theta}^{0,j}(s_{j}) \right], \quad s_{j} \in {\bm{c}}_{j}
  \end{equation}
  where $ s_{j} = \alpha_{j} \cdot \sigma_{j} = \alpha_{j} \cdot {\rm{e}}^{{\rm{i}}\theta_{j}} $, denoting the point along the mapping contour $ {\bm{c}}_{j} $ corresponding to the original point $ S_{j} $ along the physical contour $ {\bm{C}}_{j} $. The integrand of the right-hand side of Eq. (\ref{eq:2.3}) can be expanded according to Eq. (\ref{eq:2.2}) as
  \begin{equation}
    \label{eq:3.11b}
    X_{b}^{0,j}(S_{j}) + {\rm{i}} Y_{b}^{0,j}(S_{j}) = -k_{x}\gamma y(S_{j}) \frac{{\rm{d}}y(S_{j})}{{\rm{d}}S_{j}} + {\rm{i}}\gamma y(S_{j}) \frac{{\rm{d}}y(S_{j})}{{\rm{d}}S_{j}}, \quad S_{j} \in {\bm{C}}_{j}
  \end{equation}
  where $ \cos\langle \vec{n}_{j},\vec{x} \rangle = \frac{{\rm{d}}y}{{\rm{d}}S_{j}} $ and $ \cos\langle \vec{n}_{j},\vec{y} \rangle = -\frac{{\rm{d}}x}{{\rm{d}}S_{j}} $ owing to the clockwise integrating direction. The increment $ {\rm{d}}S_{j} $ in Eq. (\ref{eq:2.3}) is clockwise length increment with $ {\rm{d}}S_{j} = -|{\rm{d}}S_{j}| $, which can be mapped from the physical plane $ z = x+{\rm{i}}y $ onto mapping plane $ \zeta_{j} = \rho \cdot {\rm{e}}^{{\rm{i}}\theta_{j}} $ as
  \begin{equation}
    \label{eq:3.11c}
    {\rm{d}}S_{j} = -|{\rm{d}}S_{j}| = -|z_{j}^{\prime}(s_{j})| \cdot |{\rm{d}}s_{j}| = -|z_{j}^{\prime}(s_{j})| \cdot \alpha_{j} \cdot |\sigma_{j}| \cdot |{\rm{d}}\theta_{j}| = |z_{j}^{\prime}(s_{j})| \cdot \alpha_{j} {\rm{d}}\theta_{j}
  \end{equation}
  where $ |{\rm{d}}\theta_{j}| = - {\rm{d}}\theta_{j} $, since conformal mapping would not alter the clockwise integration direction.
\end{subequations}

With Eq. (\ref{eq:3.11a}) and (\ref{eq:3.11c}), the left-hand side of Eq. (\ref{eq:2.3}) can be rearranged as
\begin{subequations}
  \label{eq:3.12}
  \begin{equation}
    \label{eq:3.12a}
    {\rm{i}}\int_{\bm{C}_{j}} \left[ X_{a}^{0,j}(S_{j}) + {\rm{i}} Y_{a}^{0,j}(S_{j}) \right]{\rm{d}}S_{j} = \int_{\bm{c}_{j}} z_{j}^{\prime}(s_{j}) \left[ \sigma_{\rho}^{0,j}(s_{j}) + {\rm{i}} \tau_{\rho\theta}^{0,j}(s_{j}) \right] \cdot \frac{|z_{j}^{\prime}(s_{j})|}{|z_{j}^{\prime}(s_{j})|} \cdot \alpha_{j} \cdot {\rm{e}}^{{\rm{i}}\theta_{j}} \cdot {\rm{i}} \cdot {\rm{d}}\theta_{j}
  \end{equation}
  where $ {\rm{d}}s_{j} = \alpha_{j} \cdot {\rm{e}}^{\rm{i}\theta_{j}} \cdot {\rm{i}} \cdot {\rm{d}}\theta_{j} $. With Eq. (\ref{eq:3.11b}), the right-hand side of Eq. (\ref{eq:2.3}) can be rearranged as
  \begin{equation}
    \label{eq:3.12b}
    - {\rm{i}} \int_{\bm{C}_{j}} \left[ X_{b}^{0,j}(S_{j}) + {\rm{i}} Y_{b}^{0,j}(S_{j}) \right]{\rm{d}}S_{j} = - \gamma \int_{\bm{C}_{j}} y(S_{j}) \left[ k_{x} {\rm{i}}{\rm{d}}y(S_{j}) + {\rm{d}}x(S_{j}) \right]
  \end{equation}
\end{subequations}
Via the backward conformal mapping in Eq. (\ref{eq:3.1b}), the right-hand side of Eq. (\ref{eq:3.12b}) can be mapped as
\begin{equation*}
  \left\{
    \begin{aligned}
      & {\rm{d}}y(S_{j}) = \frac{{\rm{d}}y(s_{j})}{{\rm{d}}\sigma_{j}}{\rm{d}}\sigma_{j} \\
      & {\rm{d}}x(S_{j}) = \frac{{\rm{d}}x(s_{j})}{{\rm{d}}\sigma_{j}}{\rm{d}}\sigma_{j} \\
    \end{aligned}
  \right.
\end{equation*}

The righ-hand sides of Eqs. (\ref{eq:3.12a}) and (\ref{eq:3.12b}) should be equal according to Eq. (\ref{eq:2.3}), and we should have
\begin{equation}
  \label{eq:3.13}
  \int_{\bm{c}_{j}} z_{j}^{\prime}(s_{j}) \left[ \sigma_{\rho}^{0,j}(s_{j}) + {\rm{i}} \tau_{\rho\theta}^{0,j}(s_{j}) \right]{\rm{d}}s_{j} = - \gamma \int_{\bm{c}_{j}} y(s_{j}) \left[ k_{x} \frac{{\rm{i}}{\rm{d}}y(s_{j})}{{\rm{d}}\sigma_{j}} + \frac{{\rm{d}}x(s_{j})}{{\rm{d}}\sigma_{j}} \right] {\rm{d}}\sigma_{j}
\end{equation}
where
\begin{equation*}
  \left\{
    \begin{aligned}
      & x(s_{j}) = \frac{1}{2}\left[ \overline{z_{j}(s_{j})} + z_{j}(s_{j}) \right] \\
      & y(s_{j}) = \frac{{\rm{i}}}{2}\left[ \overline{z_{j}(s_{j})} - z_{j}(s_{j}) \right] \\
    \end{aligned}
  \right.
\end{equation*}
Eqs. (\ref{eq:3.10}) and (\ref{eq:3.13}) are the mixed boundary conditions mapped from Eqs. (\ref{eq:2.6}) and (\ref{eq:2.3}) into the mapping plane $ \zeta_{j} = \rho_{j} \cdot {\rm{e}}^{{\rm{i}}\theta_{j}} $, respectively. Eq. (\ref{eq:3.10}) would form a homogenerous Riemann-Hilbert problem with extra constraint in Eq. (\ref{eq:3.13}), and the solution procedure would be briefly presented below.

\section{Solution procedure of Riemann-Hilbert problem}
\label{sec:4}

Before any discussion of the solution, we should emphasize that all the discussion is conducted in the mapping plane $ \zeta_{j} = \rho_{j} \cdot {\rm{e}}^{{\rm{i}}\theta_{j}} $. The general solution of the homogenerous Riemann-Hilbert problem in Eq. (\ref{eq:3.10}) can be given below. Eq. (\ref{eq:3.10a}) would be simultaneously satisfied owing to the analytic continuation of Eq. (\ref{eq:3.6}), and no further discussion should be needed. As for Eq. (\ref{eq:3.10b}), we construct the following expression to simulate one of its potential solutions as
\begin{equation}
  \label{eq:4.1}
  X_{j}(\zeta_{j}) = (\zeta_{j}-t_{1,j})^{-\frac{1}{2}-{\rm{i}}\lambda} (\zeta_{j}-t_{2,j})^{-\frac{1}{2}+{\rm{i}}\lambda}, \quad \lambda = \frac{\ln\kappa}{2\pi}, \quad \zeta_{j} \in \overline{\bm{\omega}}_{j} \cup {\bm{\omega}}_{j}^{-}
\end{equation}
where $ t_{1,j} $ and $ t_{2,j} $ are the singularities mapped from $ T_{1} $ and $ T_{2} $ in the lower half geomaterial region $ \overline{\bm{\varOmega}}_{j} $ via the bidirectional conformal mapping in Eq. (\ref{eq:3.1}).

However, our problem consists of not only the mixed boundaries along ground surface in Eq. (\ref{eq:3.10}), but also the traction boundary condition in Eq. (\ref{eq:3.13}). The left-hand side of Eq. (\ref{eq:3.13}) contains both complex potentials $ \varphi_{0,j}^{\prime}(\zeta_{j}) $ and $ \psi_{0,j}^{\prime}(\zeta_{j}) $ to be determined, which should be analytic within the mapping region $ \overline{\bm{\omega}}_{j} $. To simulate the possible solution, we may assume the general solution of $ \varphi_{0,j}^{\prime} (\zeta_{j}) $ as a combination of $ X_{j}(\zeta_{j}) $ and an arbitrary analytic function, which should be defined within the mapping region $ \overline{\bm{\omega}}_{j} $ with two poles of origin and complex infinity. Therefore, the general solution of $ \varphi_{0,j}^{\prime} (\zeta_{j}) $ can be expressed as
\begin{equation}
  \label{eq:4.2}
  \varphi_{0,j}^{\prime} (\zeta_{j}) = X_{j}(\zeta_{j}) \sum\limits_{n = -\infty}^{\infty} f_{j,n} \zeta_{j}^{n}, \quad \zeta_{j} \in \overline{\bm{\omega}}_{j}
\end{equation}
where $ f_{j,n} $ denote complex coefficients for $ j^{\rm{th}} $ cavity to be determined. Comparing to Eq. (\ref{eq:4.1}), the definition domain of item $ X_{j}(\zeta_{}) $ in Eq. (\ref{eq:4.2}) is reduced from $ \overline{\bm{\omega}}_{j} \cup {\bm{\omega}}_{j}^{-} $ to $ \overline{\bm{\omega}}_{j} $.

The complex items with complex fractional power $ X_{j}(\zeta_{j}) $ in Eq. (\ref{eq:4.1}) is difficult to reach computational results, and can be subsequently expanded into Taylor series near the poles origin and complex infinity, respectively, as
\begin{equation}
  \label{eq:4.3}
  X_{j}(\zeta_{j}) =
  \left\{
    \begin{aligned}
      & \sum\limits_{k=0}^{\infty} \alpha_{j,k} \zeta_{j}^{k}, \quad \zeta_{j} \in {\bm{\omega}}_{j}^{+} \\
      & \sum\limits_{k=0}^{\infty} \beta_{j,k} \zeta_{j}^{-k}, \quad \zeta_{j} \in {\bm{\omega}}_{j}^{-} \\
    \end{aligned}
  \right.
\end{equation}
where
\begin{subequations}
  \label{eq:4.4}
  \begin{equation}
    \label{eq:4.4a}
    \left\{
      \begin{aligned}
        & \alpha_{j,0} = -t_{1,j}^{-\frac{1}{2}-{\rm{i}}\lambda} t_{2,j}^{-\frac{1}{2}+{\rm{i}}\lambda} \\
        & \alpha_{j,k} = -t_{1,j}^{-\frac{1}{2}-{\rm{i}}\lambda} t_{2,j}^{-\frac{1}{2}+{\rm{i}}\lambda} \cdot (-1)^{k} \left( c_{k}t_{1,j}^{-k} + \overline{c}_{k}t_{2,j}^{-k} + \sum\limits_{l=1}^{k-1} c_{l}\overline{c}_{k-l} \cdot t_{1,j}^{-l}t_{2,j}^{-k+l} \right), k \geq 1
      \end{aligned}
    \right.
  \end{equation}
  \begin{equation}
    \label{eq:4.4b}
    \left\{
      \begin{aligned}
        & \beta_{j,0} = 0 \\
        & \beta_{j,1} = 1 \\
        & \beta_{j,k} = (-1)^{k-1}\left( c_{k-1}t_{1,j}^{k-1} + \overline{c}_{k-1}t_{2,j}^{k-1} + \sum\limits_{l=1}^{k-2} c_{l}\overline{c}_{k-1-l} \cdot t_{1,j}^{l}t_{2,j}^{k-1-l} \right), \quad k \geq 2
      \end{aligned}
    \right.
  \end{equation}
  with
  \begin{equation}
    \label{eq:4.4c}
    c_{k} = \prod\limits_{l=1}^{k} \left( \frac{1}{2} - {\rm{i}}\lambda -l \right)/k!
  \end{equation}
\end{subequations}

Substituting Eq. (\ref{eq:4.3}) into Eq. (\ref{eq:4.2}) yields
\begin{subequations}
  \label{eq:4.5}
  \begin{equation}
    \label{eq:4.5a}
    \varphi_{0,j}^{\prime}(\zeta_{j}) = \sum\limits_{k=-\infty}^{\infty} A_{j,k} \zeta_{j}^{k}, \quad \zeta_{j} \in {\bm{\omega}}_{j}^{+}
  \end{equation}
  \begin{equation}
    \label{eq:4.5b}
    \varphi_{0,j}^{\prime}(\zeta_{j}) = \sum\limits_{k=-\infty}^{\infty} B_{j,k} \zeta_{j}^{k}, \quad \zeta_{j} \in {\bm{\omega}}_{j}^{-}
  \end{equation}
\end{subequations}
where
\begin{subequations}
  \label{eq:4.6}
  \begin{equation}
    \label{eq:4.6a}
    A_{j,k} = \sum\limits_{n=-\infty}^{k} \alpha_{j,k-n} f_{j,n}
  \end{equation}
  \begin{equation}
    \label{eq:4.6b}
    B_{j,k} = \sum\limits_{n=k}^{\infty} \beta_{j,-k+n} f_{j,n}
  \end{equation}
\end{subequations}
Eq. (\ref{eq:4.5}) provides expanding expressions of complex potential $ \varphi_{0,j}^{\prime}(\zeta_{j}) $ with undetermined coefficients $ f_{j,n} $ in Eq. (\ref{eq:4.2}).

Now we may provide formation of the other complex potential $ \psi_{0,j}^{\prime}(\zeta_{j}) $ using Eqs. (\ref{eq:4.5a}) and (\ref{eq:4.5b}) below. Eq. (\ref{eq:3.6}) can be formulated as
\begin{equation}
  \label{eq:4.7}
  \psi_{0,j}^{\prime}(t_{j}) = \frac{\overline{t}_{j}}{t_{j}}\overline{\varphi_{0,j}^{\prime}(t_{j})} + \frac{\overline{t}_{j}}{t_{j}}\frac{\overline{z_{j}^{\prime}(t_{j})}}{z_{j}^{\prime}(t_{j})}\varphi_{0,j}^{\prime}(t_{j}) + \frac{\overline{z_{j}(t_{j})}z_{j}^{\prime\prime}(t_{j})}{\left[ z_{j}^{\prime}(t_{j}) \right]^{2}}\varphi_{0,j}^{\prime}(t_{j}) - \frac{\overline{z_{j}(t_{j})}}{z_{j}^{\prime}(t_{j})}\varphi_{0,j}^{\prime\prime}(t_{j}), \quad t_{j} \in {\bm{c}}_{0f,j}
\end{equation}
Eq. (\ref{eq:4.7}) can be rewritten into a more compact manner as
\begin{equation}
  \label{eq:4.7'}
  \tag{4.7'}
  \psi_{0,j}^{\prime}(t_{j}) = \frac{\overline{t_{j}}}{t_{j}}\overline{\varphi_{0,j}^{\prime}(t_{j})} - {\rm{d}}\left[ \frac{\overline{z_{j}(t_{j})}}{z_{j}^{\prime}(t_{j})}\varphi_{0,j}^{\prime}(t_{j}) \right]/{\rm{d}}t_{j}, \quad t_{j} \in {\bm{c}}_{0f,j}
\end{equation}
If we note that $ \overline{z_{j}(t_{j})} = z_{j}(t_{j}) $ denotes the $ x $ axis in the physical plane $ z = x + {\rm{i}}y $ and $ \overline{t}_{j} = t_{j}^{-1} $ on the segment $ {\bm{c}}_{0f,j} $ of the unit annuli $ \overline{\bm{\omega}}_{j} $, Eq. (\ref{eq:4.7'}) can be transformed as
\begin{equation}
  \label{eq:4.8}
  \psi_{0,j}^{\prime}(t_{j}) = \frac{1}{t_{j}^{2}}\overline{\varphi}_{0,j}^{\prime}(t_{j}^{-1}) - {\rm{d}}\left[ \frac{z_{j}(t_{j})}{z_{j}^{\prime}(t_{j})}\varphi_{0,j}^{\prime}(t_{j}) \right]/{\rm{d}}t_{j}, \quad t_{j} \in {\bm{c}}_{0f,j}
\end{equation}
In above deductions, the closed-form expression of the other complex potential $ \psi_{0,j}^{\prime}(\zeta_{j}) $ is still not given, we may assume the following formation as
\begin{equation}
  \label{eq:4.9}
  \psi_{0,j}^{\prime}(\zeta_{j}) = \frac{1}{\zeta_{j}^{2}}\overline{\varphi}_{0,j}^{\prime}(\zeta_{j}^{-1}) - {\rm{d}}\left[ \frac{z_{j}(\zeta_{j})}{z_{j}^{\prime}(\zeta_{j})} \varphi_{0,j}^{\prime}(\zeta_{j}) \right]/{\rm{d}}\zeta_{j}, \quad \zeta_{j} \in \overline{\bm{\omega}}_{j}
\end{equation}
Apparently, when $ \zeta_{j} \rightarrow t_{j} \in {\bm{c}}_{0f,j} $, Eq. (\ref{eq:4.9}) would turn into its boundary value in Eq. (\ref{eq:4.8}) inversely using $ t_{j}^{-1} = \overline{t}_{j} $ in the first item and $ z_{j}(t_{j}) = \overline{z_{j}(t_{j})} $ in the second item, respectively. When $ \zeta_{j} \in \overline{\bm{\omega}}_{j} \setminus \left\{ t_{1,j},t_{2,j} \right\} $, $ \zeta_{j}^{-1} \in {\bm{\omega}}_{j}^{-} \cup {\bm{c}}_{0c,j} \cup {\bm{c}}_{0f,j} $, and the first item of Eq. (\ref{eq:4.9}) can be expanded using Eq. (\ref{eq:4.5b}) as
\begin{equation}
  \label{eq:4.9'}
  \tag{4.9'}
  \psi_{0,j}^{\prime}(\zeta_{j}) = \sum\limits_{k=-\infty}^{\infty} \overline{B}_{j,-k-2}\zeta_{j}^{k} - {\rm{d}}\left[ \frac{z_{j}(\zeta_{j})}{z_{j}^{\prime}(\zeta_{j})} \varphi_{0,j}^{\prime}(\zeta_{j}) \right]/{\rm{d}}\zeta_{j}, \quad \zeta_{j} \in \overline{\bm{\omega}}_{j}
\end{equation}
Eq. (\ref{eq:4.9'}) provides expressions of complex potential $ \psi_{0,j}^{\prime}(\zeta_{j}) $ with undetermined coefficients $ f_{j,n} $ in Eq. (\ref{eq:4.2}) in the same definition domain. In the following deductions, we should determine coefficients $ f_{j,n} $ in Eq. (\ref{eq:4.2}) using the boundary condition along cavity boundary in Eq. (\ref{eq:3.13}).

Substituting Eqs. (\ref{eq:3.3b}), (\ref{eq:4.5a}), and (\ref{eq:4.9'}) into the left-hand side of Eq. (\ref{eq:3.13}) yields
\begin{equation}
  \label{eq:4.10}
  \begin{aligned}
    & \sum\limits_{\substack{k=-\infty \\ k \neq 0}}^{\infty} \left( A_{j,k-1}\frac{\alpha_{j}^{k}}{k}\sigma_{j}^{k} + B_{j,-k-1} \frac{\alpha_{j}^{k}}{k}\sigma_{j}^{-k} \right) + \frac{z_{j}(\alpha_{j}\sigma_{j})-\overline{z_{j}(\alpha_{j}\sigma_{j})}}{\overline{z_{j}^{\prime}(\alpha_{j}\sigma_{j})}} \sum\limits_{k=-\infty}^{\infty} \overline{A}_{j,k}\alpha_{j}^{k}\sigma_{j}^{-k} \\
    & + (A_{j,-1}+B_{j,-1})\ln\alpha_{j} + C_{a,j} + (A_{j,-1}-B_{j,-1}){\rm{Ln}}\sigma_{j} = -\gamma\int_{{\bm{c}}_{j}} y(s_{j}) \left[ k_{x}\frac{{\rm{i}}{\rm{d}}y(s_{j})}{{\rm{d}}\sigma_{j}} + \frac{{\rm{d}}x(s_{j})}{{\rm{d}}\sigma_{j}} \right]{\rm{d}}\sigma_{j}
  \end{aligned}
\end{equation}
where $ C_{a,j} $ denotes an arbitrary complex integral constant, $ {\rm{Ln}} $ denotes the possibly multi-valued natural logarithm. The mapping item and the integrand of the right-hand side of Eq. (\ref{eq:3.13}) can be expanded using the sample point method \cite{lin2024over-under-excavation,lin2024charge-simulation} as
\begin{subequations}
  \label{eq:4.11}
  \begin{equation}
    \label{eq:4.11a}
    \frac{z_{j}(\alpha_{j}\sigma_{j})-\overline{z_{j}(\alpha_{j}\sigma_{j})}}{\overline{z_{j}^{\prime}(\alpha_{j}\sigma_{j})}} = \sum\limits_{k=-\infty}^{\infty} d_{j,k}\sigma_{j}^{k}
  \end{equation}
  \begin{equation}
    \label{eq:4.11b}
    y(s_{j}) \left[ k_{x}\frac{{\rm{i}}{\rm{d}}y(s_{j})}{{\rm{d}}\sigma_{j}} + \frac{{\rm{d}}x(s_{j})}{{\rm{d}}\sigma_{j}} \right] = \sum\limits_{k=-\infty}^{\infty} h_{j,k}\sigma_{j}^{k}
  \end{equation}
\end{subequations}
Then the right-hand side of Eq. (\ref{eq:3.13}) can be integrated using Eq. (\ref{eq:4.11b}) as
\begin{equation}
  \label{eq:4.12}
  -\gamma\int_{{\bm{c}}_{j}} \left[ k_{x}\frac{{\rm{i}}{\rm{d}}y(s_{j})}{{\rm{d}}\sigma_{j}} + \frac{{\rm{d}}x(s_{j})}{{\rm{d}}\sigma_{j}} \right]{\rm{d}}\sigma_{j} = \sum\limits_{\substack{k=-\infty \\ k \neq 0}}^{\infty} I_{j,k}\sigma_{j}^{k} + I_{j,0}{\rm{Ln}}\sigma_{j}
\end{equation}
where
\begin{equation}
  \label{eq:4.13}
  \left\{
    \begin{aligned}
      & I_{j,k} = -\gamma \frac{h_{j,k-1}}{k}, \quad k \neq 0 \\
      & I_{j,0} = -\gamma h_{j,-1}
    \end{aligned}
  \right.
\end{equation}

Substituting Eqs. (\ref{eq:4.10}), (\ref{eq:4.11a}), and (\ref{eq:4.12}) into Eq. (\ref{eq:3.13}) and comparing the coefficients of same power of $ \sigma_{j} $ yields
\begin{subequations}
  \label{eq:4.14}
  \begin{equation}
    \label{eq:4.14a}
    A_{j,-k-1} = -k\alpha_{j}^{k}I_{j,-k} + \alpha_{j}^{2k}B_{j,-k-1} + k\alpha_{j}^{2k} \sum\limits_{l=-\infty}^{\infty} \alpha_{j}^{l}d_{j,k}\overline{A}_{j,l+k}, \quad k \geq 1
  \end{equation}
  \begin{equation}
    \label{eq:4.14b}
    B_{j,k-1} = -k\alpha_{j}^{k}I_{j,k} + \alpha_{j}^{2k}A_{j,k-1} + k \sum\limits_{l=-\infty}^{\infty} \alpha_{j}^{l} d_{j,k} \overline{A}_{j,l-k}, \quad k \geq 1
  \end{equation}
\end{subequations}
\begin{equation}
  \label{eq:4.15}
  A_{j,-1} - B_{j,-1} = I_{j,0}
\end{equation}
The resultant equilibrium in Eq. (\ref{eq:4.15}) can be transformed using residual theorem as \cite{LIN2024appl_math_model,LIN2024comput_geotech_1,lin2024over-under-excavation,lin2024charge-simulation}
\begin{equation}
  \label{eq:4.15'}
  \tag{4.15'}
  A_{j,-1} - B_{j,-1} = -{\rm{i}}\frac{R_{y}^{0,j}}{2\pi}
\end{equation}
Eqs. (\ref{eq:4.15}) and (\ref{eq:4.15'}) reveal an implicit equilibrium as
\begin{equation}
  \label{eq:4.16}
  I_{j,0} = -\frac{\rm{i}\gamma}{2\pi} \iint_{\bm{D}_{j}}{\rm{d}}x{\rm{d}}y
\end{equation}
Eq. (\ref{eq:4.16}) should always be numerically examined. Furthermore, Eqs. (\ref{eq:4.14}) and (\ref{eq:4.15}) derived from Eqs. (\ref{eq:3.10}) and (\ref{eq:3.13}) only determines the first derivatives $ \varphi_{0,j}^{\prime}(\zeta_{j}) $ and $ \psi_{0,j}^{\prime}(\zeta_{j}) $ to be analytic and single-valued, however, the displacement in Eq. (\ref{eq:3.3c}), which contains the original functions $ \varphi_{0,j}(\zeta_{j}) $ and $ \psi_{0,j}(\zeta_{j}) $, should be also analytic and single-valued. To guarantee the single-valuedness of displacement in Eq. (\ref{eq:3.3c}), the follow equilibrium should be established as \cite{LIN2024appl_math_model,LIN2024comput_geotech_1,lin2024over-under-excavation,lin2024charge-simulation}
\begin{equation}
  \label{eq:4.17}
  \kappa A_{j,-1} + B_{j,-1} = 0
\end{equation}
Eqs. (\ref{eq:4.15}) and (\ref{eq:4.17}) give the following invariants unaffected by conformal mappings as
\begin{subequations}
  \label{eq:4.18}
  \begin{equation}
    \label{eq:4.18a}
    A_{j,-1} = \frac{I_{j,0}}{1+\kappa}
  \end{equation}
  \begin{equation}
    \label{eq:4.18b}
    B_{j,-1} = \frac{-\kappa I_{j,0}}{1+\kappa}
  \end{equation}
\end{subequations}

Substituting Eq. (\ref{eq:4.6}) into the left-hand sides of Eqs. (\ref{eq:4.14}) and (\ref{eq:4.18}) gives well-defined simultaneous complex linear system on $ f_{n} $ ($ -\infty < n < \infty $) as
\begin{subequations}
  \label{eq:4.19}
  \begin{equation}
    \label{eq:4.19a}
    \left\{
      \begin{aligned}
        & \sum\limits_{n=1}^{\infty} \alpha_{j,-1+n} f_{j,n} = \frac{I_{j,0}}{1+\kappa} \\
        & \sum\limits_{n=k+1}^{\infty} \alpha_{j,-k-1+n} f_{j,-n} = -k\alpha_{j}^{k}I_{j,-k} + \alpha_{j}^{2k}B_{j,-k-1} + k\alpha_{j}^{k} \sum\limits_{l=-\infty}^{\infty} \alpha_{j}^{l}d_{j,k}\overline{A}_{j,l+k}, \quad k \geq 1 \\
      \end{aligned}
    \right.
  \end{equation}
  \begin{equation}
    \label{eq:4.19b}
    \left\{
      \begin{aligned}
        & \sum\limits_{n=0}^{\infty} \beta_{j,1+n}f_{j,n} = \frac{-\kappa I_{j,0}}{1+\kappa} \\
        & \sum\limits_{n=k}^{\infty} \beta_{j,-k+1+n}f_{j,n} = 	-k\alpha_{j}^{k}I_{j,k} + \alpha_{j}^{2k}A_{j,k-1} + k \sum\limits_{l=-\infty}^{\infty} \alpha_{j}^{l} d_{j,k} \overline{A}_{j,l-k}, \quad k \geq 1 \\
      \end{aligned}
    \right.
  \end{equation}
\end{subequations}
Eq. (\ref{eq:4.19}) can be solved in the iterative method in Refs \cite{LIN2024appl_math_model,LIN2024comput_geotech_1,lin2024over-under-excavation,lin2024charge-simulation} to uniquely determine $ f_{n} $ in Eq. (\ref{eq:4.2}), so that the complex potentials in Eqs. (\ref{eq:4.5a}) and (\ref{eq:4.9'}) can be subsequently determined for further solution of stress and displacement of sequential shallow tunnelling.

\section{Stress and displacement fields of sequential shallow tunnelling}
\label{sec:5}

Replacing $ s_{j} = \alpha_{j}\sigma_{j} $ with $ \zeta_{j} = \rho_{j}\sigma_{j} $ in Eqs. (\ref{eq:3.13}) and (\ref{eq:4.10}) yields the integration of the curvilinear traction components mapped onto the mapping annuli $ \overline{\bm{\omega}}_{j}\setminus\left\{ t_{1,j},t_{2,j} \right\} $ for a given polar radius $ \rho_{j} $ as
\begin{equation}
  \label{eq:5.1}
  \begin{aligned}
    & \int z_{j}^{\prime}(\rho_{j}\sigma_{j})\left[ \sigma_{\rho}^{0,j}(\rho_{j}\sigma_{j}) + {\rm{i}} \tau_{\rho\theta}^{0,j}(\rho_{j}\sigma_{j}) \right]{\rm{d}}(\rho_{j}\sigma_{j}) \\
    = 
    & \sum\limits_{\substack{k=-\infty \\ k \neq 0}}^{\infty} \left( A_{j,k-1}\frac{\rho_{j}^{k}}{k} - B_{j,k-1}\frac{\rho_{j}^{-k}}{k} \right) \sigma_{j}^{k} + \sum\limits_{k=-\infty}^{\infty} \sum\limits_{l=-\infty}^{\infty} g_{j,l}(\rho_{j})\overline{A}_{j,l-k}\rho_{j}^{l-k}\sigma_{j}^{k} \\
    & + \left( A_{j,-1} + B_{j,-1} \right)\ln\rho_{j} + \left( A_{j,-1} - B_{j,-1} \right){\rm{Ln}}\sigma_{j} + C_{a,j}, \quad \alpha_{j} \leq \rho_{j} \leq 1
  \end{aligned}
\end{equation}
where
\begin{equation}
  \label{eq:5.2}
  \sum\limits_{k=-\infty}^{\infty} g_{j,k}(\rho_{j})\sigma_{j}^{k} = \frac{z_{j}(\rho_{j}\sigma_{j})-\overline{z_{j}(\rho_{j}\sigma_{j})}}{\overline{z_{j}^{\prime}(\rho_{j}\sigma_{j})}}
\end{equation}
When $ \rho_{j} \rightarrow \alpha_{j} $, Eq. (\ref{eq:5.1}) would be degenerated into Eq. (\ref{eq:4.10}). When $ \rho_{j} \rightarrow 1 $, the right-hand side of Eq. (\ref{eq:5.2}) would be zero, since $ \overline{z_{j}(\sigma_{j})} = z_{j}(\sigma_{j}) $ denotes the $ x $ axis in the physical plane $ z = x + {\rm{i}}y $.

The derivation of Eq. (\ref{eq:5.1}) gives the expansion of Eq. (\ref{eq:3.3b}) as
\begin{subequations}
  \label{eq:5.3}
  \begin{equation}
    \label{eq:5.3a}
    \sigma_{\rho}^{0,j}(\rho_{j}\sigma_{j}) + {\rm{i}} \tau_{\rho\theta}^{0,j}(\rho_{j}\sigma_{j}) = \frac{1}{z_{j}^{\prime}(\rho_{j}\sigma_{j})} \sum\limits_{k=-\infty}^{\infty} \left( 
      \begin{aligned}
        & A_{j,k}\rho_{j}^{k} - B_{j,k}\rho_{j}^{-k-2} \\
        + & (k+1)\sum\limits_{l=-\infty}^{\infty} g_{j,l}(\rho_{j})\overline{A}_{j,l-k-1}\rho_{j}^{l-k-2} \\
      \end{aligned}
    \right)\sigma_{j}^{k}, \quad \alpha_{j} \leq \rho_{j} \leq 1
  \end{equation}
  Eqs. (\ref{eq:3.3a}) and (\ref{eq:3.3c}) can be respectively expanded as
  \begin{equation}
    \label{eq:3.38b}
    \sigma_{\theta}^{0,j}(\rho_{j}\sigma_{j}) + \sigma_{\rho}^{0,j}(\rho_{j}\sigma_{j}) = 4\Re\left[ \frac{1}{z_{j}^{\prime}(\rho_{j}\sigma_{j})} \sum\limits_{k=-\infty}^{\infty} A_{j,k}\rho_{j}^{k}\sigma_{j}^{k} \right], \quad \alpha_{j} \leq \rho_{j} \leq 1
  \end{equation}
  \begin{equation}
    \label{eq:3.38c}
    \begin{aligned}
      2G[u_{0,j}(\rho_{j}\sigma_{j})+{\rm{i}}v_{0,j}(\rho_{j}\sigma_{j})] = 
      & \sum\limits_{\substack{k=-\infty \\ k \neq 0}}^{\infty} \frac{1}{k} \left( \kappa A_{j,k-1}\rho_{j}^{k}\sigma_{j}^{k} + B_{j,k-1}\rho_{j}^{-k}\sigma_{j}^{k} \right) + C_{0} + (\kappa A_{j,-1} - B_{j,-1})\ln\rho_{j} \\
      & - \sum\limits_{k=-\infty}^{\infty} \sum\limits_{l=-\infty}^{\infty} g_{j,l}(\rho_{j})\overline{A}_{j,l-k}\rho_{j}^{l-k}\sigma_{j}^{k}, \quad \alpha_{j} \leq \rho_{j} \leq 1
    \end{aligned}
  \end{equation}
  with
  \begin{equation*}
    C_{0} = -\sum\limits_{\substack{k=-\infty \\ k \neq 0}}^{k=\infty}\frac{1}{k}\left( \kappa A_{j,k-1} + B_{j,k-1} \right)
  \end{equation*}
\end{subequations}

The initial stress field in Eq. (\ref{eq:2.1}) can be mapped onto the mapping annuli $ \overline{\bm{\omega}}_{j} $ as
\begin{equation}
  \label{eq:2.1'}
  \tag{2.1'}
  \left\{
    \begin{aligned}
      & \sigma_{\theta}^{0}(\rho_{j}\sigma_{j}) + \sigma_{\rho}^{0}(\rho_{j}\sigma_{j}) = \sigma_{y}^{0}[z_{j}(\rho_{j}\sigma_{j})] + \sigma_{x}^{0}[z_{j}(\rho_{j}\sigma_{j})] \\
      & \sigma_{\theta}^{0}(\rho_{j}\sigma_{j}) - \sigma_{\rho}^{0}(\rho_{j}\sigma_{j}) + 2{\rm{i}}\tau_{\rho\theta}^{0}(\rho_{j}\sigma_{j}) = \left\{ \sigma_{y}^{0}[z_{j}(\rho_{j}\sigma_{j})] - \sigma_{x}^{0}[z_{j}(\rho_{j}\sigma_{j})] + 2{\rm{i}}\tau_{xy}[z(\rho_{j}\sigma_{j})] \right\} \cdot \frac{z_{j}^{\prime}(\rho_{j}\sigma_{j})}{\overline{z_{j}^{\prime}(\rho_{j}\sigma_{j})}} \cdot \sigma_{j}^{2}
    \end{aligned}
  \right., \quad \alpha_{j} \leq \rho_{j} \leq 1
\end{equation}
where $ \sigma_{\theta}^{0} $, $ \sigma_{\rho}^{0} $, and $ \tau_{\rho\theta}^{0} $ denote hoop, radial, and shear components of the initial stress field mapped onto the mapping unit annuli $ \overline{\bm{\omega}}_{j} $, respectively.

According to Eqs. (\ref{eq:5.3}) and (\ref{eq:2.1'}), the total curvilinear stress field mapped onto mapping unit annuli $ \overline{\bm{\omega}}_{j} $ can be obtained as
\begin{subequations}
  \label{eq:5.5}
  \begin{equation}
    \label{eq:5.5a}
    \left\{
      \begin{aligned}
        & \sigma_{\theta}^{j}(\rho_{j}\sigma_{j}) = \sigma_{\theta}^{0}(\rho_{j}\sigma_{j}) + \sigma_{\theta}^{0,j}(\rho_{j}\sigma_{j}) \\
        & \sigma_{\rho}^{j}(\rho_{j}\sigma_{j}) = \sigma_{\rho}^{0}(\rho_{j}\sigma_{j}) + \sigma_{\rho}^{0,j}(\rho_{j}\sigma_{j}) \\
        & \tau_{\rho\theta}^{j}(\rho_{j}\sigma_{j}) = \tau_{\rho\theta}^{0}(\rho_{j}\sigma_{j}) + \tau_{\rho\theta}^{0,j}(\rho_{j}\sigma_{j}) \\
      \end{aligned}
    \right., \quad \alpha_{j} \leq \rho_{j} \leq 1
  \end{equation}
  where $ \sigma_{\theta}^{j} $, $ \sigma_{\rho}^{j} $, and $ \tau_{\rho\theta}^{j} $ denote hoop, radial, and tangential components of total curvilinear stress mapped onto the mapping unit annuli $ \overline{\bm{\omega}}_{j} $ along the circles of radii $ |\zeta_{j}| = \rho_{j} $, respectively. The displacement field mapped onto mapping unit annuli $ \overline{\bm{\omega}}_{j} $ can be obtained as
  \begin{equation}
    \label{eq:5.5b}
    u_{j}(\rho_{j}\sigma_{j}) + {\rm{i}}v_{j}(\rho_{j}\sigma_{j}) = u_{j}[z_{j}(\rho_{j}\sigma_{j})] + {\rm{i}}v_{j}[z_{j}(\rho_{j}\sigma_{j})], \quad \alpha_{j} \leq \rho_{j} \leq 1, \sigma_{j} \neq t_{1,j},t_{2,j}
  \end{equation}
\end{subequations}

Significantly, when $ \rho_{j} \rightarrow \alpha_{j} $, Eq. (\ref{eq:5.5a}) gives the Mises stress (numerically equal to the absolute value of $ \sigma_{\theta}^{j}(\alpha_{j}\sigma_{j}) $) and the residual stresses (radial component $ \sigma_{\rho}^{j}(\alpha_{j}\sigma_{j}) $ and tangential component $ \tau_{\rho\theta}^{j}(\alpha_{j}\sigma_{j}) $) along cavity boundary $ {\bm{C}}_{j} $ caused by $ j^{\rm{th}} $ stage excavation; Eq. (\ref{eq:5.5b}) gives the horizontal and vertical deformation along cavity boundary $ {\bm{C}}_{j} $.

The curvilinear stress and displacement fields in Eq. (\ref{eq:5.5}) can be mapped onto rectangular ones as
\begin{subequations}
  \label{eq:5.6}
  \begin{equation}
    \label{eq:5.6a}
    \left\{
      \begin{aligned}
        & \sigma_{y}^{j}[z_{j}(\zeta_{j})] + \sigma_{x}^{j}[z_{j}(\zeta_{j})] = \sigma_{\theta}^{j}(\zeta_{j}) + \sigma_{\rho}^{j}(\zeta_{j}) \\
        & \sigma_{y}^{j}[z_{j}(\zeta_{j})] - \sigma_{x}^{j}[z_{j}(\zeta_{j})] + 2{\rm{i}}\tau_{xy}^{j}[z_{j}(\zeta_{j})] = \left[ \sigma_{\theta}^{j}(\zeta_{j}) - \sigma_{\rho}^{j}(\zeta_{j}) + 2{\rm{i}}\tau_{\rho\theta}^{j}(\zeta_{j}) \right] \cdot \frac{\overline{\zeta}_{j}}{\zeta_{j}} \cdot \frac{\overline{z_{j}^{\prime}(\zeta_{j})}}{z_{j}^{\prime}(\zeta_{j})}
      \end{aligned}
    \right.
    , \quad \zeta_{j} = \rho_{j} \cdot \sigma_{j}
  \end{equation}
  \begin{equation}
    \label{eq:5.6b}
    u_{j}[z_{j}(\zeta_{j})] + {\rm{i}}v_{j}[z_{j}(\zeta_{j})] = u_{j}(\zeta_{j}) + v_{j}(\zeta_{j}), \quad \zeta_{j} = \rho_{j}\cdot\sigma_{j}
  \end{equation}
\end{subequations}
where $ \sigma_{x}^{j} $, $ \sigma_{y}^{j} $, and $ \tau_{xy}^{j} $ denote horizontal, vertical, and shear components of the total stress field after $ j^{\rm{th}} $-stage excavation, respectively.

\section{Numerical case and verification}
\label{sec:6}

The solution in Section \ref{sec:5} is analytical with infinite Laurent series of $ f_{j,n} $ in Eq. (\ref{eq:4.2}), which should be bilaterally truncated into $ 2M_{j}+1 $ items to reach numerical results. The coefficient series in Eq. (\ref{eq:4.5}) should be truncated as
\begin{equation}
  \label{eq:4.5a'}
  \tag{4.5a'}
  A_{j,k} = \sum\limits_{n=-M_{j}}^{k} \alpha_{j,k+n} f_{j,-n}
\end{equation}
\begin{equation}
  \label{eq:4.5b'}
  \tag{4.5b'}
  B_{j,k} = \sum\limits_{n=k}^{M_{j}} \beta_{j,-k+n} f_{j,n}
\end{equation}
where $ j = 1,2,3,\cdots,N $. The infinite series to obtain solution of $ f_{j,n} $ in Eqs. (\ref{eq:4.2}) should be truncated as well. Correspondingly, the stress and displacement fields in Section \ref{sec:5} should be truncated, which leads to numerical errors in oscillation form of Gibbs phenomena \cite{LIN2024appl_math_model,lin2024over-under-excavation}. To reduce the Gibbs phenomena, the Lanczos filtering is applied in Eq. (\ref{eq:5.3}) to replace $ A_{j,k} $ and $ B_{j,k} $ with $ A_{j,k} \cdot L_{k} $ and $ B_{j,k} \cdot L_{k} $, respectively, where
\begin{equation}
  \label{eq:6.1}
  L_{k} = \left\{
    \begin{aligned}
      & 1, \quad k = 0 \\
      & \sin\left( \frac{k}{M_{j}}\pi \right)/\left( \frac{k}{M_{j}}\pi \right), \quad {\rm{otherwise}}
    \end{aligned}
  \right.
\end{equation}
with $ -M_{j} \leq k \leq M_{j} $.

In this section, a numerical case of 4 stage shallow tunnelling is performed to validate the present analytical solution by comparing with a corresponding finite element solution. The present analytical solution is conducted using the programming language \texttt{Fortran} of compiler \texttt{GCC} (version14.1.1). The linear systems are solved using the \texttt{dgesv} and \texttt{zgesv} packages of \texttt{LAPACK} (version 3.12.0). The figures are plotted using \texttt{gnuplot} (version 6.0).

\subsection{Cavity geometries and bidirectional conformal mappings of sequential shallow tunnelling}
\label{sec:6.1}

The numerical case of sequential shallow tunnelling takes the geometry of the 4-stage excavation in Fig. \ref{fig:2}b, while the sharp corner smoothening technique is applied to adapt the numerical schemes of bidirectional conformal mapping in \ref{sec:a2}. Thus, the cavity boundaries of 4-stage excavation can be analytically expressed as
\begin{subequations}
  \label{eq:6.2}
  \begin{equation}
    \label{eq:6.2a}
    {\rm{Stage \; 1}}: \quad
    \left\{
      \begin{aligned}
        & x^{2} + (y+10)^{2} = 5^{2}, \quad -5 \leq x \leq 0, -10 \leq y \leq -5 \\
        & (x+4.5)^{2} + (y+10)^{2} = 0.5^{2}, \quad -5 \leq x \leq -4.5, -10.5 \leq y \leq -10 \\
        & y = -10.5, \quad -4.5 \leq x \leq 0 \\
        & x^{2} + (y+10)^{2} = 0.5^{2}, \quad 0 \leq x \leq 0.5, -10.5 \leq y \leq -10 \\
        & x = 0.5, \quad -10 \leq y \leq -5.5 \\
        & x^{2} + (y+5.5)^{2} = 0.5^{2}, \quad 0 \leq x \leq 0.5, -5.5 \leq y \leq -5 \\
      \end{aligned}
    \right.
  \end{equation}
  \begin{equation}
    \label{eq:6.2b}
    {\rm{Stage \; 2}}: \quad
    \left\{
      \begin{aligned}
        & x^{2} + (y+10)^{2} = 5^{2}, \quad -5 \leq x \leq 5, -10 \leq y \leq -5 \\
        & (x+4.5)^{2} + (y+10)^{2} = 0.5^{2}, \quad -5 \leq x \leq -4.5, -10.5 \leq y \leq -10 \\
        & y = -10.5, \quad -4.5 \leq x \leq 4.5 \\
        & (x-4.5)^{2} + (y+10)^{2} = 0.5^{2}, \quad 4.5 \leq x \leq 5, -10.5 \leq y \leq -10 \\
      \end{aligned}
    \right.
  \end{equation}
  \begin{equation}
    \label{eq:6.2c}
    {\rm{Stage \; 3}}: \quad
    \left\{
      \begin{aligned}
        & x^{2} + (y+10)^{2} = 5^{2}, \quad -5 \leq x \leq 5, -10 \leq y \leq -5 \\
        & x = -5, \quad -14.5 \leq y \leq -10 \\
        & (x+4.5)^{2} + (y+14.5)^{2} = 0.5^{2}, \quad -5 \leq x \leq -4.5, -15 \leq y \leq -14.5 \\
        & y = -15, \quad -4.5 \leq x \leq 0 \\
        & x^{2} + (y+14.5)^{2} = 0.5^{2}, \quad 0 \leq x \leq 0.5, -15 \leq y \leq -14.5 \\
        & x = 0.5, \quad -14.5 \leq y \leq -11 \\
        & (x-1)^{2} + (y+11)^{2} = 0.5^{2}, \quad 0.5 \leq x \leq 1, -11 \leq y \leq -10.5 \\
        & y = -10.5, \quad 1 \leq x \leq 4.5 \\
        & (x-4.5)^{2} + (y+10)^{2} = 0.5^{2}, \quad 4.5 \leq x \leq 5, -10.5 \leq y \leq -10 \\
      \end{aligned}
    \right.
  \end{equation}
  \begin{equation}
    \label{eq:6.2d}
    {\rm{Stage \; 4}}: \quad
    \left\{
      \begin{aligned}
        & x^{2} + (y+10)^{2} = 5^{2}, \quad -5 \leq x \leq 5, -10 \leq y \leq -5 \\
        & x = -5, \quad -14.5 \leq y \leq -10 \\
        & (x+4.5)^{2} + (y+14.5)^{2} = 0.5^{2}, \quad -5 \leq x \leq -4.5, -15 \leq y \leq -14.5 \\
        & y = -15, \quad -4.5 \leq x \leq 4.5 \\
        & (x-4.5)^{2} + (y+14.5)^{2} = 0.5^{2}, \quad 4.5 \leq x \leq 5, -15 \leq y \leq -14.5 \\
        & x = 5, \quad -14.5 \leq y \leq -10 \\ 
      \end{aligned}
    \right.
  \end{equation}
\end{subequations}

\begin{figure}[htb]
  \centering
  \includegraphics[width = 0.8\textwidth]{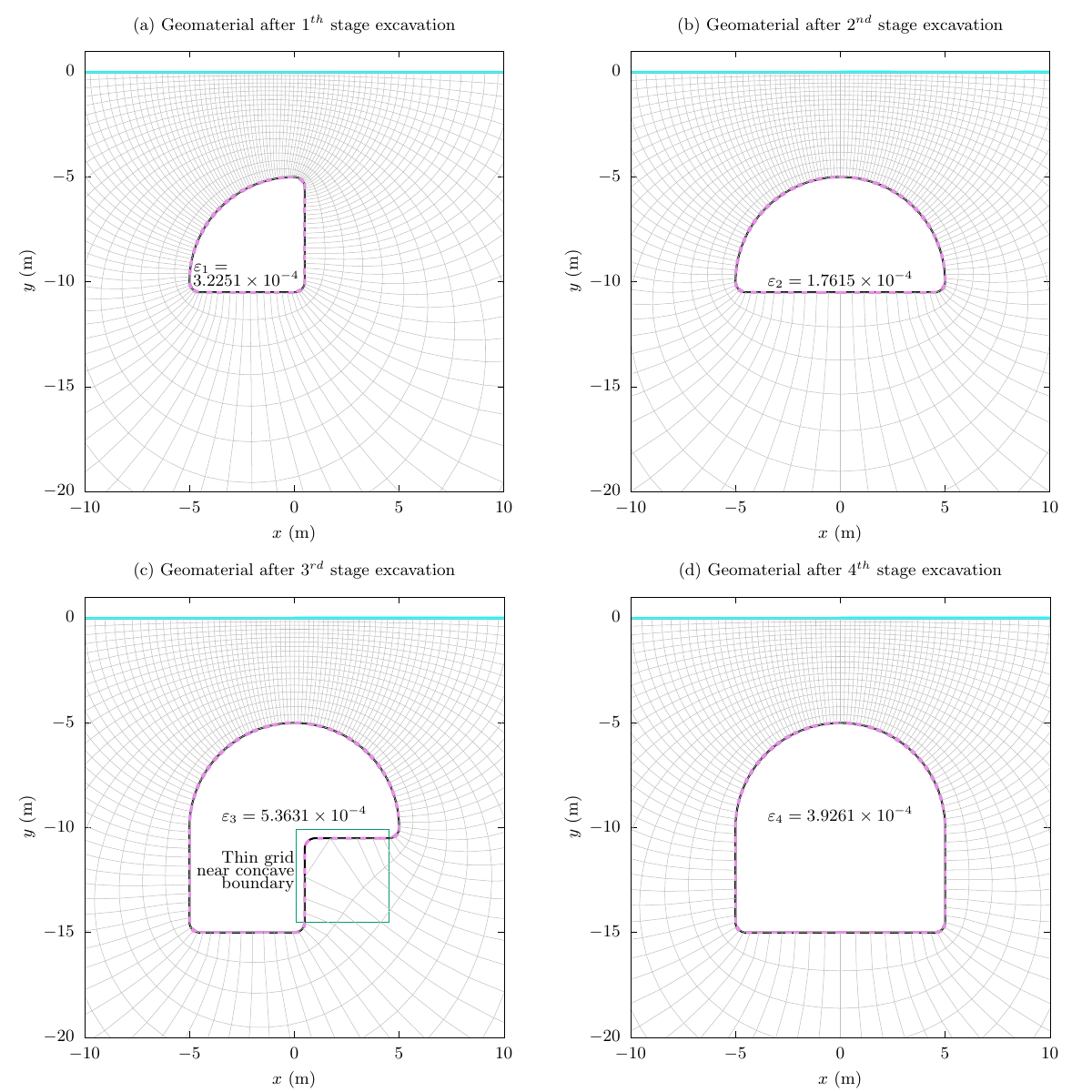}
  \caption{Curvilinear grids of backward conformal mapping of geomaterial for four excavation stages}
  \label{fig:4}
\end{figure}

With the geometrical expressions in Eq. (\ref{eq:6.2}), the bidirectional conformal mappings of the 4-stage sequential shallow tunnelling can be obtained in Fig. \ref{fig:4}, in which the collocation points are taken in the following pattern for good accuracy after trial computations: (1) 30 points are uniformly selected along the $ 90^{\circ} $ arc of $ 0.5{\rm{m}} $ radius; (2) 90 points are uniformly selected along the $ 90^{\circ} $ arc of $ 5{\rm{m}} $ radius; 60 points are uniformly selected along the rest straight lines. The collocation points along the cavity boundary in the physical plane $ z $ are denoted by $ z_{j,i}^{f} $, where $ j $ denotes the excavation stage, $ i $ denotes the collocation point number, $ f $ denotes the word "forward". After trial computations, the assignment factors take $ K_{j}^{0} = k_{j}^{0} = 2.2 $ and $ K_{j}^{1} = k_{j}^{1} = 1.5 $, while $ w_{j}^{c} = -2.5-7.5{\rm{i}} $ and $ w_{j,0} = 0.9\beta $, where $ \beta = 5 $ and $ j = 1,2,3,4 $.

To quantitatively determine the accuracy of the bidirectional conformal mapping, the following critieria are established as
\begin{equation}
  \label{eq:6.3}
  \varepsilon_{j} = \max | z_{j,i}^{m,b} - z_{j,i}^{m,f} |, \quad i = 1,2,3,\cdots,N_{j}
\end{equation}
where $ z_{j,i}^{m,f} $ denote the midpoints of the selected collocation points $ z_{j,i}^{f} $, and can be computed as
\begin{equation*}
  z_{j,i}^{m,f} = \frac{1}{2}(z_{j,i}^{f}+z_{j,i+1}^{f}), \quad i = 1,2,3,\cdots,N_{j}
\end{equation*}
$ z_{j,i}^{m,b} $ denote the coordinates of the corresponding midpoints $ z_{j,i}^{m,f} $ after forward and backward conformal mappings, which can be computed in the following procedure.

Via Eq. (\ref{eqa:1a}), the midpoints of the collocation points $ z_{j,i}^{m,f} $ can be forwardly mapped onto corresponding mapping points $ w_{j,i}^{m,f} $ in the interval mapping plane $ \bm{\varOmega}_{j}^{w} $. Then via Eq. (\ref{eqa:2a'}), the midpoints in the mapping plane $ \bm{\omega}_{j} $ can be computed as $ \zeta_{j,i}^{m,f} $. Subsequently, the midpoints for backward conformal mapping in the mapping plane $ \bm{\omega}_{j} $ are given by only preserving the polar angles as
\begin{equation}
  \label{eq:6.4}
  \zeta_{j,i}^{m,b} = \alpha_{j} \cdot \exp({\rm{i}}\arg\zeta_{j,i}^{m,f})
\end{equation}
where $ b $ in the superscript denotes "backward". Apparently, slight differences between $ \zeta_{j,i}^{m,b} $ and  $ \zeta_{j,i}^{m,f} $ should exist. Whereafter, the backward midpoints in the interval mapping plane $ \bm{\varOmega}_{j}^{w} $ can be computed as $ w_{j,i}^{m,b} $ via Eq. (\ref{eqa:6}). Finally, the backward midpoints $ z_{j,i}^{m,b} $ in the physical plane $ \bm{\varOmega}_{j} $ can be computed via Eq. (\ref{eqa:1b}). Therefore, $ \varepsilon_{j} $ in Eq. (\ref{eq:4.3}) should be able to present the maximum difference between the original midpoints and computed ones along cavity boundary for sequential excavation stages.

The orthogonal grids and corresponding values of $ \varepsilon_{j} $ in Fig. \ref{fig:4} show that the bidirectional conformal mapping in \ref{sec:a} is very adaptive, and can be used to conformally map a lower half plane containing an asymmetrical cavity of very complicated shape (see Fig. \ref{fig:4}c). Therefore, the conformal mapping can be used in further mechanical computation. However, we should also note in Fig. \ref{fig:4}c that the orthogonal grid near concave boundary would be thin (as marked), indicating that the accuracy for a lower half plane containing a concave cavity would be slightly compromised.

\subsection{Comparisons with finite element solution}
\label{sec:6.2}

\begin{figure}[htb]
  \centering
  \includegraphics[width = \textwidth]{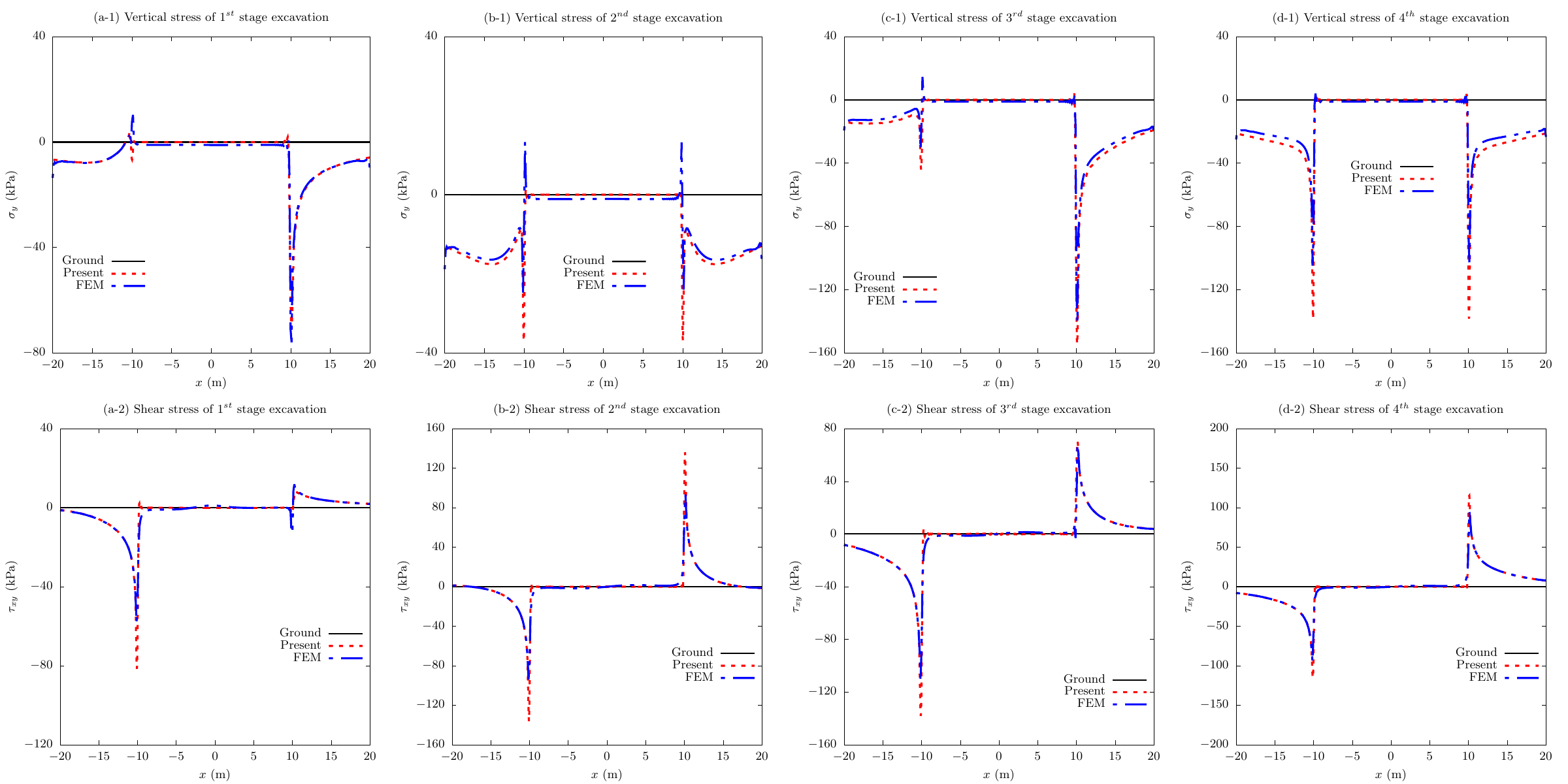}
  \caption{Comparisons of vertical and shear stress components along ground surface of four excavation stages between present solution and finite element solution}
  \label{fig:8}
\end{figure}

\begin{figure}[htb]
  \centering
  \includegraphics[width = \textwidth]{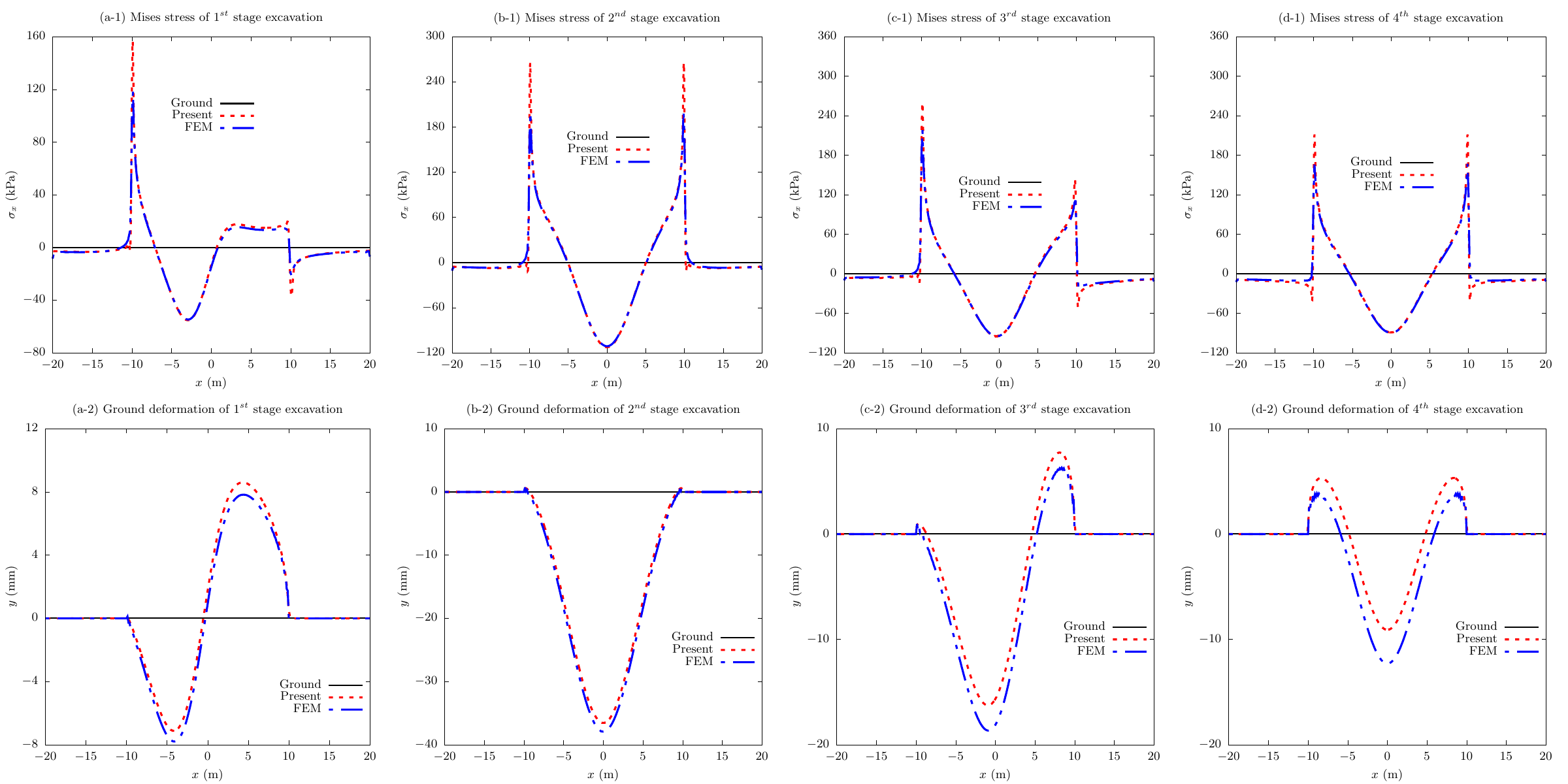}
  \caption{Comparisons of horizontal stress and deformation along ground surface of four excavation stages between present solution and finite element solution}
  \label{fig:9}
\end{figure}

To validate the present analytical solution in Section \ref{sec:5}, a finite element solution using software \texttt{ABAQUS 2020} is performed for comparisons. The mechanical parameters of geomaterial are listed in Table \ref{tab:1}, where the coordinates of the joint points $ T_{1} $ and $ T_{2} $ are deliberately located near the shallow tunnel for better illustration of boundary conditions along ground surface (see Figs. \ref{fig:8} and \ref{fig:9}).

\begin{table}[htb]
  \centering
  \scriptsize
  \caption{Mechanical parameters of numerical case}
  \label{tab:1}
  \begin{tabular}{cccccc}
    \toprule
    $ \gamma {\rm{(kN/m^{3})}} $  & $ k_{x} $ & $ E {\rm{(MPa)}} $ & $ \nu $ & $ x_{0} $ (m) & $ M_{j} $ \\
    \midrule
    $ 20 $ & $ 0.8 $ & $ 20 $ & $ 0.3 $ & $ 10 $ & $ 250 $ \\
    \bottomrule
  \end{tabular}
\end{table}

\begin{figure}[htb]
  \centering
  \includegraphics[width = 0.8\textwidth]{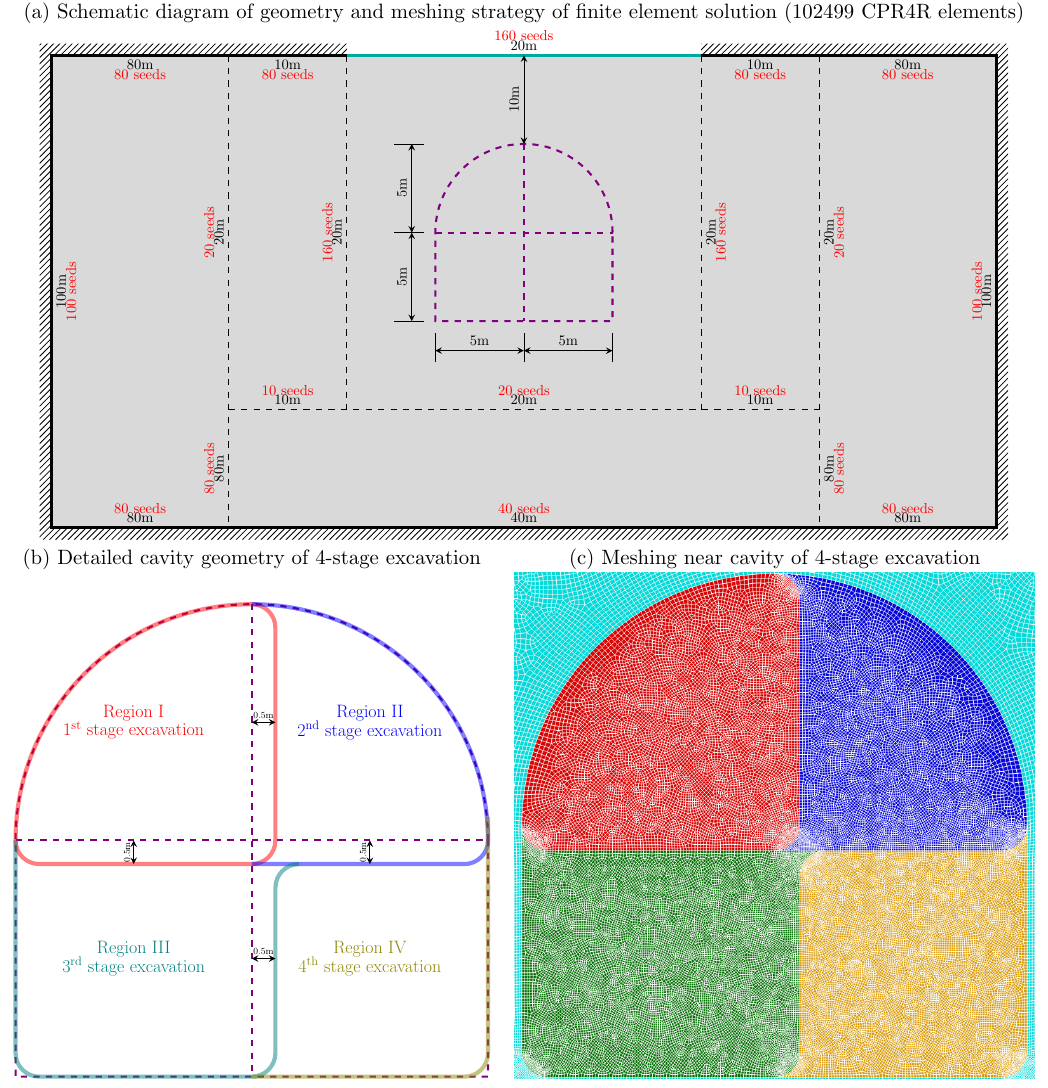}
  \caption{Geometry and meshing scheme of the finite element solution}
  \label{fig:5}
\end{figure}

The $ {\rm{kN-m}} $ mechanical model of the finite element solution is shown in Fig. \ref{fig:5}a, where the same 4-stage sequential shallow tunnelling is conducted. In Fig. \ref{fig:5}a, the model geometry, geomaterial sketching, and seed distribution scheme outside of cavity are illustrated. In Fig. \ref{fig:5}b, the 4-stage excavation procedure is marked by different colors, and more detailed seed distribution scheme of cavity boundaries is elaborated below: (1) 90 seeds are uniformly distributed along the $ 90^{\circ} $ arcs of $ 5{\rm{m}} $; (2) 30 seeds are uniformly distributed along the $ 90^{\circ} $ arcs of $ 0.5{\rm{m}} $; (3) 80 seeds are uniformly distributed along the rest straight lines. The meshing near cavity is shown in Fig. \ref{fig:5}c, where the 4-stage excavation regions are marked with the same colors with the diagram of Fig. \ref{fig:5}b. Both analytical and finite element solutions are in plane strain condition, thus, the 102499 finite elements use CPR4R type. The steps of finite element solution are listed in Table \ref{tab:2}, where the 4 excavation stages are sequentially conducted, according to Fig. \ref{fig:5}b.

\begin{table}[htb]
  \centering
  \scriptsize
  \caption{Steps of finite element solution}
  \label{tab:2}
  \begin{tabular}{ccccccc}
    \toprule
    Step & Initial & Step 0 & Step 1 & Step 2 & Step 3 & Step 4 \\
    \midrule
    Procedure & (Initial) & Geostatic & Static,General & Static,General & Static,General & Static,General \\
    Load & \makecell[c]{Applying constraints \\ and geostress} & Applying gravity & - & - & - & - \\
    Interaction & - & - & \makecell[c]{Deactivating \\ region I} & \makecell[c]{Deactivating \\ region II} & \makecell[c]{Deactivating \\ region III} & \makecell[c]{Deactivating \\ region IV} \\
    \bottomrule
  \end{tabular}
\end{table}

\begin{figure}[htb]
  \centering
  \includegraphics[width = \textwidth]{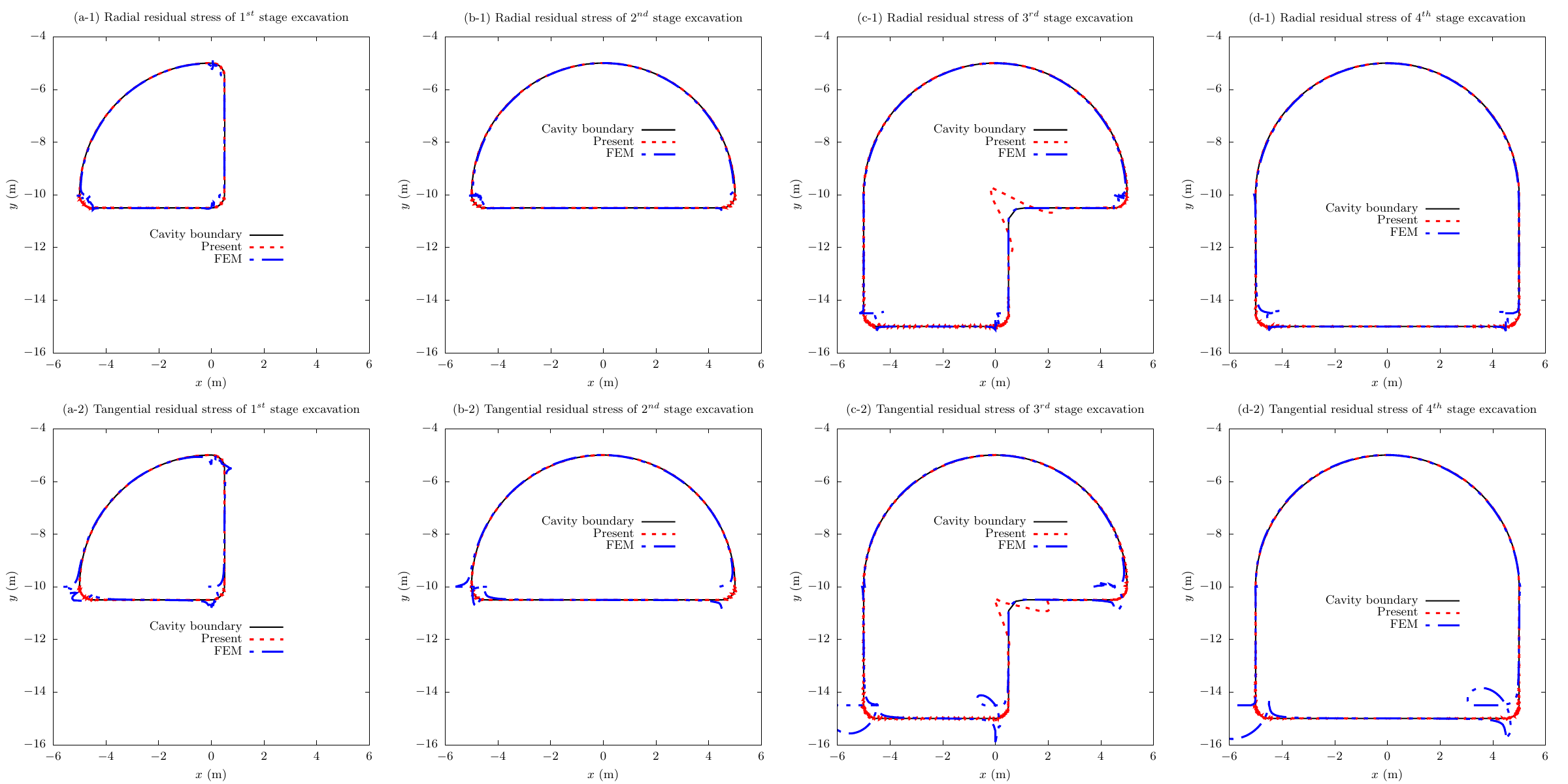}
  \caption{Comparisons of residual stress components along cavity boundaries of four excavation stages between present solution and finite element solution}
  \label{fig:6}
\end{figure}

\begin{figure}[htb]
  \centering
  \includegraphics[width = \textwidth]{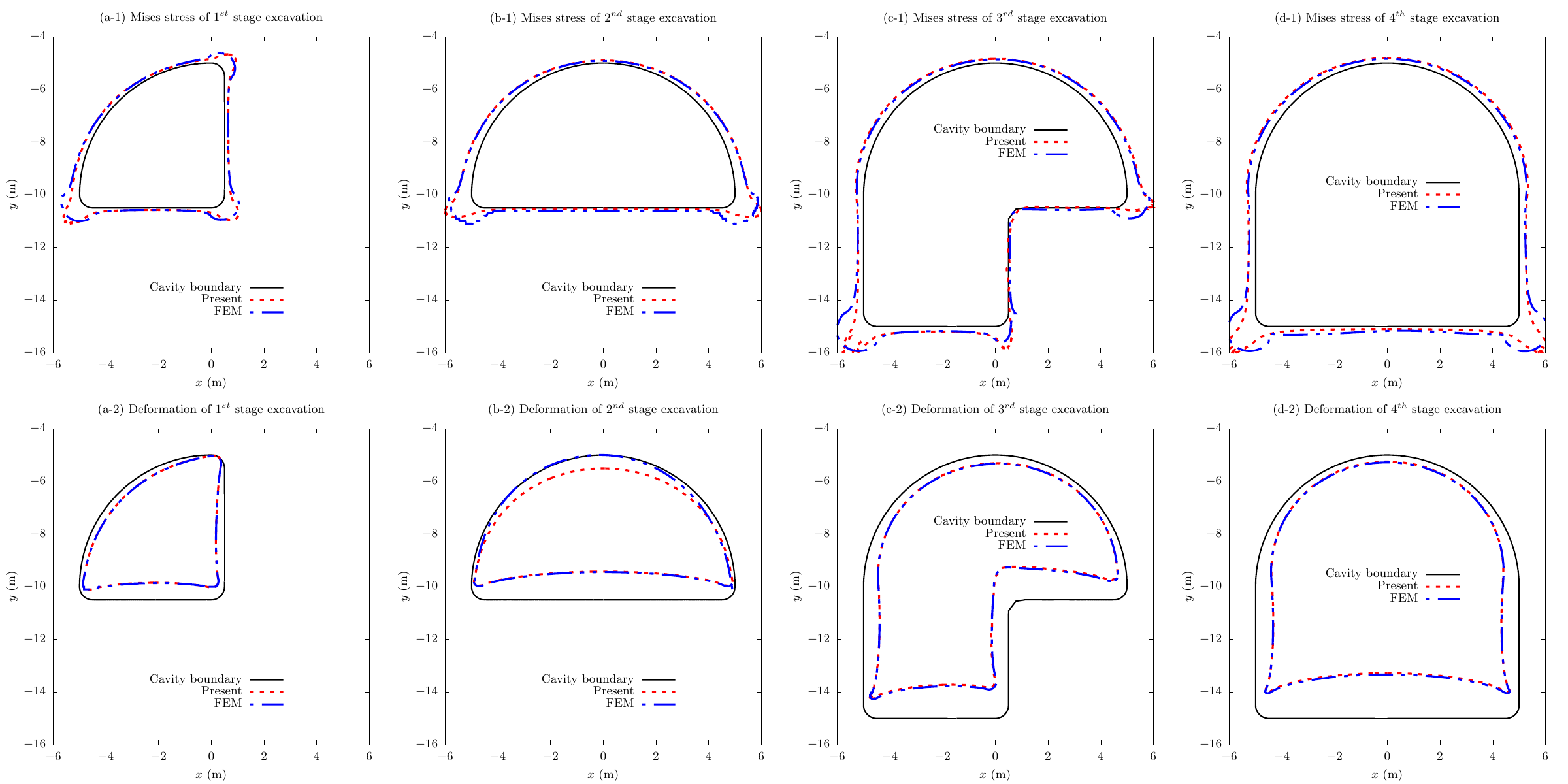}
  \caption{Comparisons of Mises stress and deformation along cavity boundaries of four excavation stages between present solution and finite elemetn solution}
  \label{fig:7}
\end{figure}

Substituting the mechanical parameters in Table \ref{tab:1} into the present solution and the finite element solution gives the results in Figs. \ref{fig:6}-\ref{fig:9}. Fig. \ref{fig:6} shows the radial and tangential components of the residual stress along cavity boundaries of all four excavation stages, which are reduced by $ 10^{-2} $ for better illustration. In theory, the residual stress should be zero to accurately meet the zero boundary condition along cavity boundary (see Fig. \ref{fig:3}). Fig. \ref{fig:6} indicates that both the present solution and the finite element solution agree to each other well for cavity boundaries with large curvature radii (a straight line has an infinite curvature radius). As for the rounded corners of small curvature radii, the present solution clearly surpasses the finite element solution. However, the residual stresses of $ 3^{\rm{rd}} $ stage excavation computed by the present solution contain significant errors near the concave boundary, indicating the defect of the present solution on concave cavity. Such a defect is caused by the bidirectional conformal mapping, and the curvilinear grid in Fig. \ref{fig:4}c might provide certain clues. In Fig. \ref{fig:4}c, the curvilinear grid density near the concave boundary is sparse and thin, which possibly causes loss of nearby mathematical information on the mechanical solution based on the curvilinear grid. Therefore, the present solution is more accurate than the finite element method for convex cavities.

Fig. \ref{fig:7} shows the Mises stress (reduced by $ 10^{-3} $) and deformation (magnified by $ 10^{2} $) along cavity boundaries for four excavation stages, and good agreements between these two solution are observed, except that the Mises stresses near the rounded corners obtained by the present solution are acuter and larger than those obtained by the finite element solution. The results of the present solution should be more accurate than the finite element solution, since the residual stresses near the rounded corners obtained by the finite element in Fig. \ref{fig:6} contains much larger errors.

Fig. \ref{fig:8} shows the comparisons of vertical and shear stress components along ground surface between the present solution and the finite element solution in all four excavation stages, and good agreements are observed. To be specific, the zero tractions in the range $ x \in [-10,10] $ meet the requirement of the boundary condition in Eq. (\ref{eq:2.6b}), further indicating the correctness of the present solution.

Fig. \ref{fig:9} shows the comparisons of horizontal stress and ground deformation between the present solution and the finite element solution. The zero deformation in the range $ x \in (\infty,-10] \cup [10,\infty) $ meets the requirement of the boundary condition in Eq. (\ref{eq:2.6a}), indicating that the present solution is capable of obtaining fixed far-field displacement as desired. The horizontal stress in the range $ x \in (\infty,-10] \cup [10,\infty) $ between the present solution and the finite element solution is in good agreements. However, the horizontal stress and ground deformation in the range $ x \in [-10,10] $ between present solution and finite element solution show clear discrepancy, especially for $ 1^{\rm{st}} $, $ 2^{\rm{nd}} $, and $ 3^{\rm{rd}} $ excavation stages. Noting that the overall indices $ \varepsilon_{j} (j = 1,2,3,4) $ are small with order of magnitutde $ 10^{-4} $ in Fig. \ref{fig:4}, the bidirectional conformal mapping should be accurate, and such stress and ground deformation discrepancies are not caused by the mapping. The residual radial and tangential stress components in Fig. \ref{fig:6} (except for \ref{fig:6}c), the vertical and shear stress components in the range $ [-10,10] $ in Fig. \ref{fig:8}, and the ground deformation in the range $ (-\infty,-10] \cup [10,\infty) $ in Fig. \ref{fig:9} match the corresponding mixed boundary conditions in Eqs. (\ref{eq:2.3}), (\ref{eq:2.6b}), and (\ref{eq:2.6a}), respectively. By contrast, the residual radial and tangential stress components of the finite element solution in Fig. \ref{fig:6} contain significant errors, especially around the rounded corners. Therefore, it would be reasonable to presume that the results computed by the present solution should be more accurate than those computed by the finite element method. In other words, the present solution should be more accurate than the finite element solution. However, we should also notice that the obvious errors of residual radial and tangential stress components computed by the present solution in Fig. \ref{fig:6}c, which are probably caused by the thin orthogonal grid near the concave. Such errors indicate that the present solution may be inaccurate for excavation with concave cavity.

In summary, Figs. \ref{fig:6}-\ref{fig:9} show the validation of the boundary conditions of the present solution, and the good agreements between the present solution and corresponding finite element solution. Additionally, detailed discussions reveal that the present solution should be more accurate than the finite element solution. The present solution can be used for further application.

\section{Parametric investigation}
\label{sec:7}

In this section, several practical applications of present solution are illustrated to show its possible usage in sequential shallow tunnelling. To be consistent and simple, the benchmark geometry and sequential excavation stages take the same ones in Section \ref{sec:4}(see Figs. \ref{fig:4} and \ref{fig:5}b), and the benchmark parameters take the ones in Table \ref{tab:1}.

\subsection{Deformation along cavity boundary and solution convergence}
\label{sec:7.1}

Deformation along cavity boundary is significant in shallow tunnelling to estimate possible displacement around tunnel for further design optimization. In the validation of present solution in Section \ref{sec:5}, the free ground surface is kept within the range of $ x \in [-x_{0},x_{0}] $ of $ x_{0} = 10{\rm{m}} $ for good visualization and comparisons with finite element solution. However, the free ground surface should range wider in practical shallow tunnelling, thus, we select the joint coordinate as $ x_{0} = 10^{1}, 10^{2}, 10^{3}, 10^{4} {\rm{m}} $ to investigate the possible deformation along cavity boundary, while the rest parameters remain the same to the ones in Table \ref{tab:1}.

\begin{figure}[htb]
  \centering
  \includegraphics[width = \textwidth]{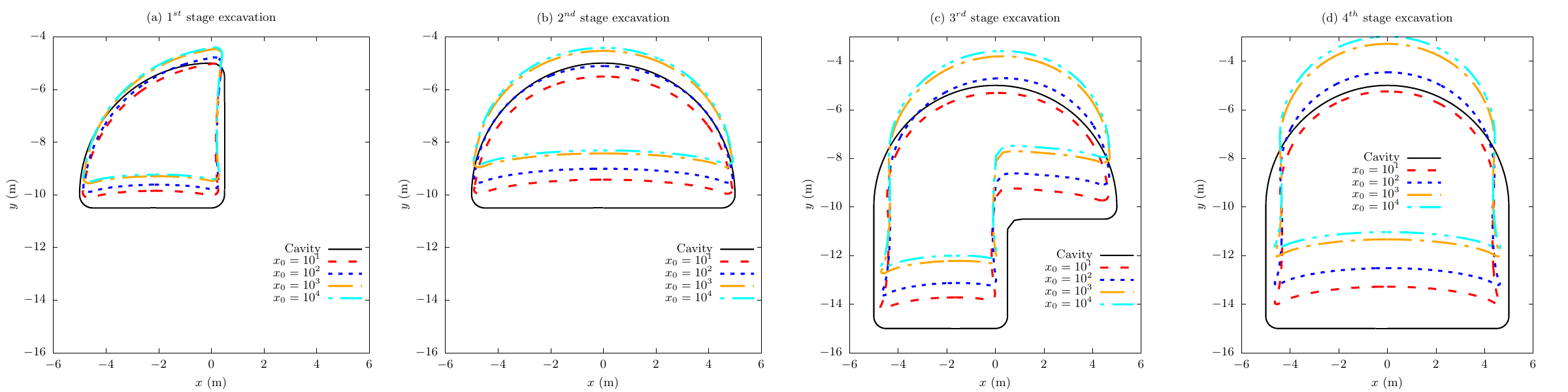}
  \caption{Cavity deformation (magnified by $ 10^{2} $ times) against joint point coordinate $ x_{0} $ for four excavation stages}
  \label{fig:10}
\end{figure}

The results of deformation along cavity boundary against joint point coordinate $ x_{0} $ for four excavation stages are shown in Fig. \ref{fig:10} with a magnification of $ 10^{2} $ times for better illustration. Fig. \ref{fig:10} shows overall upheavals of geomaterial around tunnels and clear inward contractions against cavity boundaries for all four excavation stages, which are coincident to engineering expectations. Moreover, the deformed cavity boundaries of $ x_{0} = 10^{3} {\rm{m}} $ and $ x_{0} = 10^{4} {\rm{m}} $ overlap to each other, which indicates the same deformation for large free ground segment $ {\bm{C}}_{0f} $ above tunnel, as well as convergence of present solution.

\subsection{Mises stress along cavity boundary}
\label{sec:7.2}

Mises stress along cavity boundary can be mechanically interpreted as the hoop stress along cavity boundary, which is significant in shallow tunnelling to estimate possible damage zones around cavity for consideration of necessary reinforcement measures. In the validation, the lateral coefficient is set to be $ k_{x} = 0.8 $, however, the shallow strata to excavate are complicated with different lateral coefficients. Thus, we select lateral coefficient as $ k_{x} = 0.6, 0.8, 1.0, 1.2, 1.4, 1.6 $ to investigate the variation of Mises stress along cavity boundaries of sequential shallow tunnelling, and the rest parameters remain the same to the ones in Table \ref{tab:1}.

\begin{figure}[htb]
  \centering
  \includegraphics[width = \textwidth]{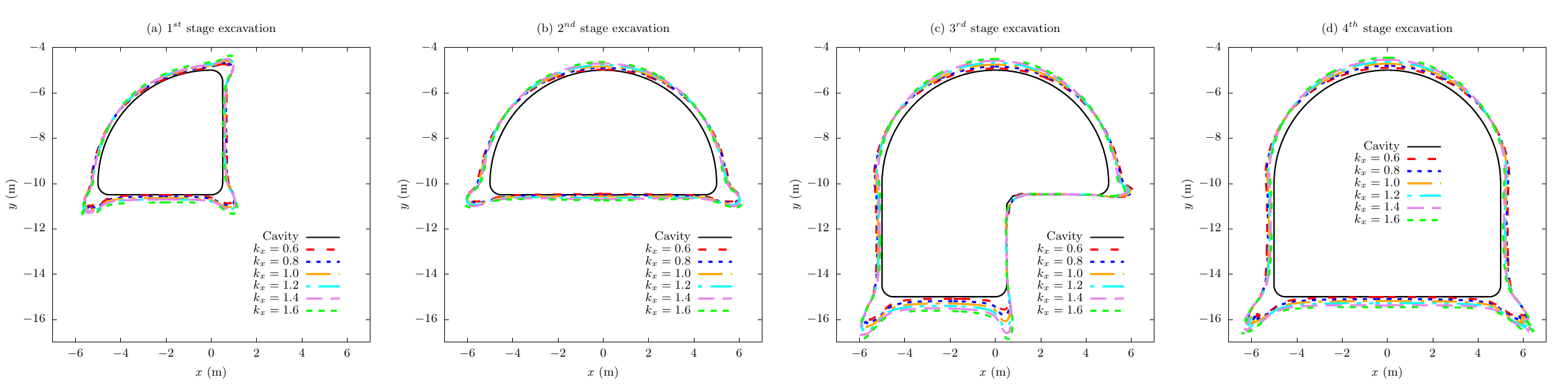}
  \caption{Mises stress along cavity boundary (reduced by $ 10^{-3} $) against lateral coeffcient $ k_{x} $ for four excavation stages}
  \label{fig:11}
\end{figure}

The results of Mises stress along cavity boundaries against lateral coefficient $ k_{x} $ for for excavation stages are shown in Fig. \ref{fig:11} with reduction of $ 10^{-3} $ times for better illustration. Fig. \ref{fig:11} shows that a larger lateral coefficient would cause higher stress concentrations near the rounded corners during sequential excavation, since the overall stress of the initial stress field increases. Therefore, tunnel corners may need more reinforcement measures in sequential excavation to avoid possible damage of geomaterial due to sress concentration. Moreover, Fig. \ref{fig:11}c shows that the Mises stress near the right bottom geomaterial is very approximate to zero, indicating possible negative hoop stress. In other words, possible tensile stress may exist near the to-be-excavated region of $ 3^{\rm{rd}} $-stage excavation, which may be hazardous in construction safety, and necessary monitoring of nearby geomaterial should be conducted.

\subsection{Stress concentration of rounded corners}
\label{sec:7.3}

Fig. \ref{fig:11} shows remarkable stress concentrations near the corners during sequential shallow tunnelling, thus, the reduction of stress concentration near corners should be discussed due to its high value in shallow tunnelling safety. Among the four excavation stages in Fig. \ref{fig:11}, the geometry of the $ 3^{\rm{th}} $ excavation stage is the most complicated, and worthy of discussion. The corners are the geometrical cause of stress concentration in Fig. \ref{fig:11}, thus, we select different radii of the rounded corners (see the excavation scheme in Fig. \ref{fig:5}b) as $ 0.3, 0.4, 0.5, 0.6, 0.7, 0.8 {\rm{m}} $, and the rest cavity geometry is the same to that in Fig. \ref{fig:4}c. The mechanical parameters take the same values in Table \ref{tab:1}.

\begin{figure}[htb]
  \centering
  \includegraphics[width = 0.75\textwidth]{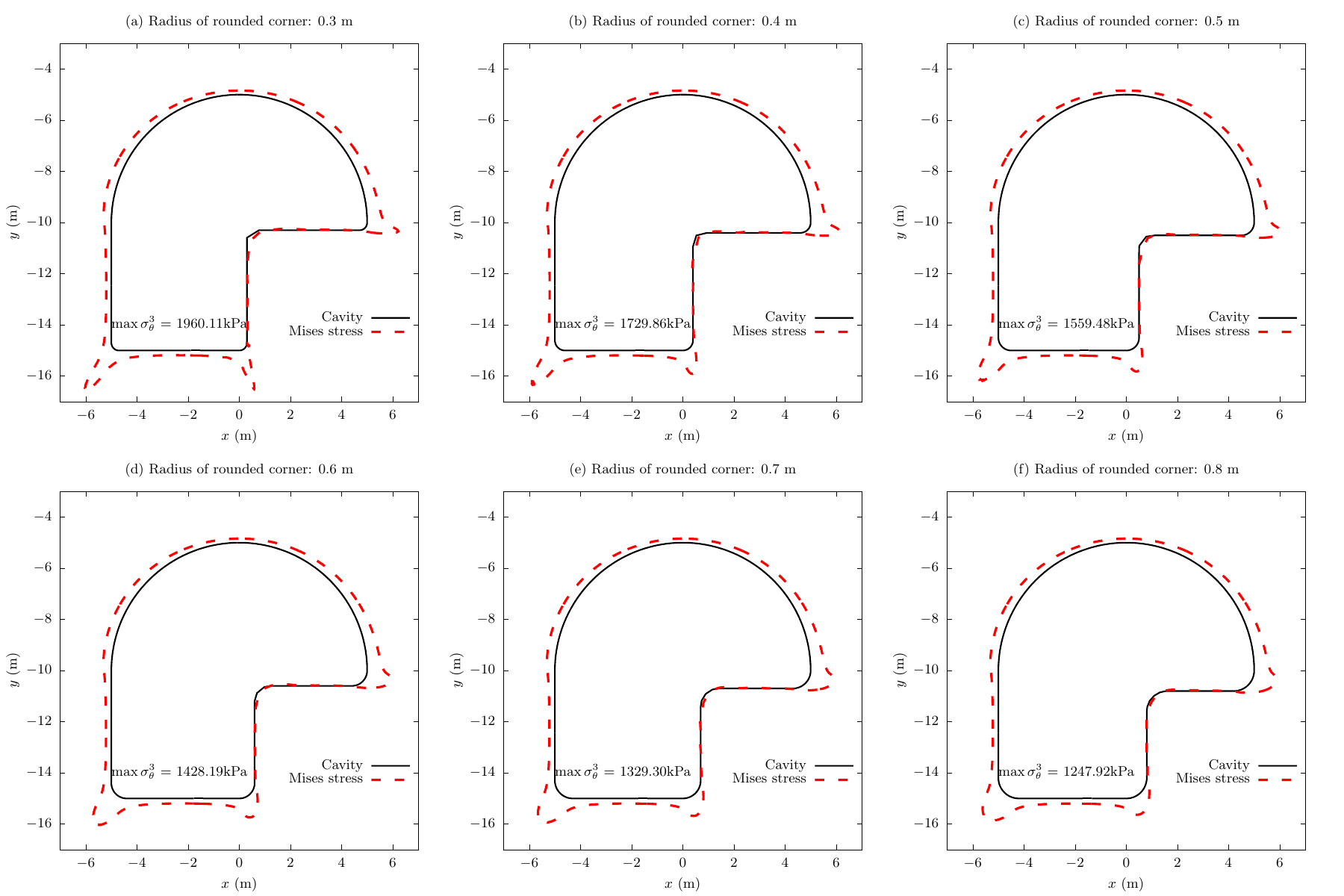}
  \caption{Mises stress along cavity boundary (reduced by $ 10^{-3} $) of different radii of rounded corner for the $ 3^{\rm{th}} $-stage excavation (maximum Mises stresses located at the left bottom corner)}
  \label{fig:12}
\end{figure}

The Mises stresses along cavity boundary of different radii of rounded corners for the $ 3^{\rm{rd}} $-stage excavation are shown in Fig. \ref{fig:12} with reduction of $ 10^{-3} $ times for better illustration. Fig. \ref{fig:12} shows that as the radius of rounded corners increases from $ 0.3{\rm{m}} $ to $ 0.8{\rm{m}} $, the stress concentration decreases from $ 1960.11{\rm{kPa}} $ to $ 1247.92{\rm{kPa}} $ by a remarkable $ 36.33\% $. Therefore, relatively large radius of rounded corners might be considered during sequential excavation to avoid overlarge stress concentration.

\section{Further discussions}
\label{sec:8}

In exsiting complex variable solutions of shallow tunnelling in gravitational geomaterial, the cavity shapes are always fully circular \cite{Lu2016,Lu2019new_solution} or symmetrically noncircular \cite{Zengguisen2019,zhou2024analytical,fan2024analytical}, and no sequential shallow tunnelling is considered. The most fundamental reason is the mathematical difficulty to seek a suitable conformal mapping scheme that could bidirectionally map a lower half plane containing a noncircular and possibly asymmetrical cavity with shape changing due to sequential excavation. The conformal mapping scheme in \ref{sec:a} extends the bidirectional stepwise conformal mapping in our previous study \cite{lin2024charge-simulation}, which incorporates Charge Simulation Method \cite{amano1994charge} and Complex Dipole Simulation Method \cite{sakakibara2020bidirectional}, to a lower half plane containing more complicated cavity shapes, as shown in Fig. \ref{fig:4} (especially in Fig. \ref{fig:4}c). Furthermore, the validation in Figs. \ref{fig:6}-\ref{fig:9} indicate the possibilities to extend the mathematical usage of the complex variable method from the full shallow tunnelling schemes with circular or symmetrically noncircular cavity shapes in the exsiting complex variable solutions \cite{Lu2016,Lu2019new_solution,Zengguisen2019,zhou2024analytical,fan2024analytical} to sequential ones with complicated interval cavity shapes (see Figs. \ref{fig:6}c and \ref{fig:7}c).

Despite the theoretical improvements mentioned above, several defects of the present solution exist, and should be disclosed as:

(1) The most obvious defect of the present solution is that the present solution requires the rest geomaterial after sequential excavation to always remain doubly-connected to hold topological consistence of geomaterial and ensure solvability of the solution, while temporary supports are not considered for seeking mechanical variation within gravitational geomaterial alone. Temporary supports are important mechanical structure in sequential shallow tunnelling, and would alter the stress and displacement fields within rest geomaterial. Ignoring temporary supports would compromise the present solution to a certain extent. In our future studies, the complex variable method should be further improved and modified to be able to consider liners and temporary supports.

(1$^{\prime}$) On the other hand, however, the numerical cases in Sections \ref{sec:4} and \ref{sec:5}, which simulate the multi-stepwise upper half vertical subdivision method in Fig. \ref{fig:2}c without consideration of temporary supports, show that the present solution is capable of dealing with complicated interval cavity shapes of sequential shallow tunnelling. Then the present solution can thereby be degenerated to adapt more simple cavity shapes, for example, the top heading and bench method in Fig. \ref{fig:2}b, where temporary supports may be not necessary. Therefore, the present solution is suitable for sequential shallow tunnelling methods that need no temporary supports without solution modification. The numerical cases in Sections \ref{sec:4} and \ref{sec:5}, which certainly deviate from tunnel engineering fact without consideration of temporary supports, are deliberately conducted to demonstrate the significant theoretical improvements of present solution.

(2) Sequential excavation generally occurs with three-dimensional effect of tunnel face, which is not considered in the present solution. Thus, the stress and displacement fields computed by the present solution would deviate from real-world stress and displacement fields to a certain extent. To be conservative and to ensure construction safety, the present solution can be only used as very preliminary estimation of sequential shallow tunnelling.

(3) The thin orthogonal grid in Fig. \ref{fig:4}c and residual stress errors in Fig. \ref{fig:6}c indicate a possible weakness of the present solution that the stress and displacement for sequential excavation stage of shallow tunnels with concave cavity may be inaccurate. For computation and design of concave cavity, the present solution may be not suitable.

\section{Conclusions}
\label{sec:9}

This paper proposes a new complex variable solution on asymmetrical sequential excavation of large-span shallow tunnels in gravitational geomaterial with consideration of rigid static equilibrium. The asymmetrical cavities of sequential tunnelling are conformally mapped using a new bidirectional stepwise mapping scheme consiting of Charge Simulation Method and Complex Dipole Simulation Method. The sequential shallow tunnelling procedure is mathematically disassembled into parallel and indepedent shallow tunnelling problems to seek similar pattern of mixed boundary conditions and solution procedure, which are subsequently transformed into a homogenerous Riemann-Hilbert problem to obtain reasonable stress and displacement field. Via comparisons with corresponding finite element solution, the bidirectional stepwise conformal mapping scheme and the proposed solution are both sufficiently validated. A final parametric investigation shows several possible applications of the proposed solution in asymmetrical sequential shallow tunnelling, and certain engineering recommendations are proposed. This new solution adapts lower half geomaterial containing asymmetrical cavity of complicated shape, and extends the usage of complex variable method in shallow tunnelling.

\section*{Acknowlegements}

This study is financially supported by the Fujian Provincial Natural Science Foundation of China (Grant No. 2022J05190), the Scientific Research Foundation of Fujian University of Technology (Grant No. GY-H-22050), and the National Natural Science Foundation of China (Grant No. 52178318). The authors would like to thank Ph.D. Yiqun Huang for his suggestion on this study.

\section*{Author Contributions}

Conceptualization: Luo-bin Lin, Fu-quan Chen; Methodology: Luo-bin Lin; Validation: Luo-bin Lin; Writing - original draft preparation: Luo-bin Lin; Writing - review and editing: Fu-quan Chen, Chang-jie Zheng, Shan-shun Lin; Funding acquisition: Chang-jie Zheng, Shang-shun Lin; Resources: Fu-quan Chen, Chang-jie Zheng; Supervision: Fu-quan Chen, Chang-jie Zheng, Shang-shun Lin

\section*{Declarations}

The authors have no relevant financial or non-financial interests to disclose.

\section*{Data Availability}

The research codes can be found via link \underline{github.com/luobinlin987/asymmetrical-sequential-static-equilibrium}.

\appendix
\section{Bidirectional conformal mapping}
\label{sec:a}

The geomaterials $ \overline{\bm{\varOmega}}_{j} $ ($ j = 1,2,3,\cdots,n $) after excavation are lower half planes containing corresponding cavities of arbitrary shapes $ {\bm{D}}_{j} $, which are doubly-connected regions and can be bidirectionally mapped onto unit annuli $ \overline{\bm{\omega}}_{j} $ using the two-step mapping scheme in our previous study \cite{lin2024charge-simulation}, respectively. Since the theory of such bidirectinal conformal mapping has been elaborated in detail in Ref \cite{lin2024charge-simulation}, the mapping scheme below only illustrates necessary numerical procedure for conciseness. Comparing to existing mapping schemes \cite{exadaktylos2002closed,Zengguisen2019,ye2023novel-aizhiyong}, the proposed bidirectional conformal mapping has the following advantages: (1) Good adaptivity for a lower half plane containing an irregular and asymmetrical cavity; (2) Fast computation of solving a pair of linear systems below; (3) Closed logic to provide both forward and backward conformal mappings for mathematical rigor.

\subsection{Mapping scheme}
\label{sec:a1}

The first step is to bidirectionally map the geomaterial $ \overline{\bm{\varOmega}}_{j} $ in the physical plane $ z = x+{\rm{i}}y $ onto interval annuli $ \overline{\bm{\varOmega}}_{j}^{w} $ in interval mapping plane $ w_{j} = \Re w_{j} + {\rm{i}} \Im w_{j} $ using the following mapping schemes as
\begin{subequations}
  \label{eqa:1}
  \begin{equation}
    \label{eqa:1a}
    w_{j} = w_{j}(z) = \beta_{j} \frac{z-z_{c}^{j}}{z-\overline{z}_{c}^{j}}, \quad z \in \overline{\bm{\varOmega}}_{j}
  \end{equation}
  \begin{equation}
    \label{eqa:1b}
    z = z_{j}(w_{j}) = \frac{w_{j}\overline{z}_{c}^{j}-\beta_{j}z_{c}^{j}}{w_{j}-\beta_{j}}, \quad w_{j} \in \overline{\bm{\varOmega}}_{j}^{w}
  \end{equation}
\end{subequations}
where $ \beta_{j} $ denotes the exterior radius of the interval mapping annulus $ \overline{\bm{\varOmega}}_{j}^{w} $, which should be given before-hand.

Via Eq. (\ref{eqa:1}), the finite free segment $ {\bm{C}}_{0f,j} $ and far-field segment $ {\bm{C}}_{0c,j} $ of the ground surface $ {\bm{C}}_{0,j} $ in the physical plane $ z = x+{\rm{i}}y $ are bidirectionally mapped onto free arc segment $ {\bm{c}}_{0f,j}^{w} $ and constrained arc segment $ {\bm{c}}_{0c,j}^{w} $ of the exterior boundary $ {\bm{c}}_{0,j}^{w} = {\bm{c}}_{0f,j}^{w} \cup {\bm{c}}_{0c,j}^{w} $ in the interval mapping plane $ w_{j} = \Re w_{j} + {\rm{i}} \Im w_{j} $, respectively. The joint points $ T_{1} $ and $ T_{2} $ connecting segments $ {\bm{C}}_{0f,j} $ and $ {\bm{C}}_{0c,j} $ in the physical plane $ z = x+{\rm{i}}y $ are bidirectionally mapped onto corresponding joint points $ t_{1,j}^{w} $ and $ t_{2,j}^{w} $ connecting arc segments $ {\bm{c}}_{0f,j}^{w} $ and $ {\bm{c}}_{0c,j}^{w} $ in the interval mapping plane $ w_{j} = \Re w_{j} + {\rm{i}} \Im w_{j} $, respectively. The cavity boundary $ {\bm{C}}_{j} $ in the physical plane $ z = x+{\rm{i}}y $ is bidirectionally mapped onto $ {\bm{c}}_{j}^{w} $ in the interval mapping plane $ w_{j} = \Re w_{j} + {\rm{i}} \Im w_{j} $. Figs. \ref{fig:a1}a and \ref{fig:a1}b graphically and schematically show the bidirectional conformal mapping in Eq. (\ref{eqa:1}) for the example in Fig. \ref{fig:3}b.

\setcounter{figure}{0}
\renewcommand{\thefigure}{A-\arabic{figure}}
\begin{figure}[htb]
  \centering
  \includegraphics[width = \textwidth]{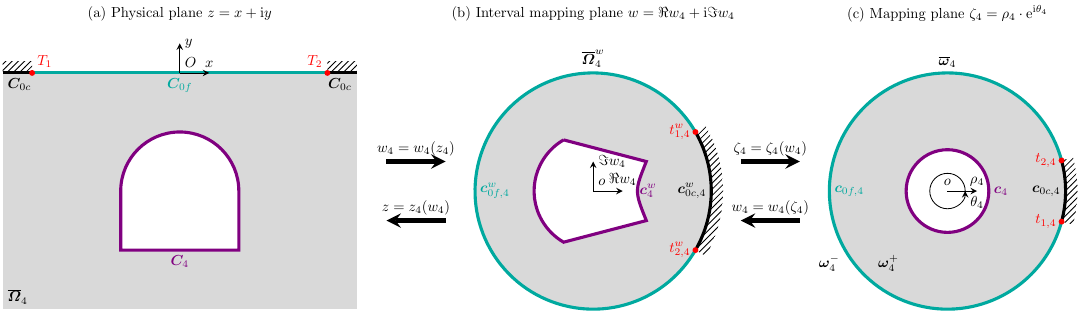}
  \caption{Stepwise bidirectional conformal mapping -- a schematic example for Fig. \ref{fig:3}b}
  \label{fig:a1}
\end{figure}

The second step is to bidirectionally map the interval annulus $ \overline{\bm{\varOmega}}_{j}^{w} $ in the interval mapping plane $ w_{j} = \Re w_{j} + {\rm{i}} \Im w_{j} $ onto corresponding unit mapping annulus $ \overline{\bm{\omega}}_{j} $ of interior radius $ \alpha_{j} $ in the mapping plane $ \zeta_{j} = \rho_{j} \cdot {\rm{e}}^{{\rm{i}}\theta_{j}} $. The forward and backward conformal mappings of second-step bidirectional conformal mapping can be respectively expressed using Charge Simulation Method \cite{amano1994charge,okano2003numerical,lin2024charge-simulation} and Complex Dipole Simulation Method \cite{sakakibara2020bidirectional,lin2024charge-simulation} as
\begin{subequations}
  \label{eqa:2}
  \begin{equation}
    \label{eqa:2a}
    \zeta_{j} = \zeta_{j}(w_{j}) = \frac{w_{j}-w_{j}^{c}}{w_{j}^{\beta}-w_{j}^{c}} \cdot \exp\left[ \sum\limits_{k=1}^{N_{0}} P_{j,k} \ln\frac{w_{j}-U_{j,k}}{w_{j}^{\beta}-U_{j,k}} + \sum\limits_{k=1}^{N_{j}} Q_{j,k} \ln\frac{w_{j}-V_{j,k}}{w_{j}^{\beta}-V_{j,k}} \right]
  \end{equation}
  \begin{equation}
    \label{eqa:2b}
    w_{j} = w_{j}(\zeta_{j}) = \sum\limits_{k=1}^{N_{0}} \frac{p_{j,k}}{\zeta_{j}-\eta_{j,k}} + \sum\limits_{k=1}^{N_{j}} \frac{q_{j,k}}{\zeta_{j}-\mu_{j,k}}
  \end{equation}
\end{subequations}
where $ w_{j}^{c} $ denotes arbitrary point within interior boundary $ {\bm{c}}_{j}^{w} $ of the interval mapping annulus $ \overline{\bm{\varOmega}}_{j}^{w} $ for fixed point normalization; $ w_{j}^{\beta} $ denotes arbitrary point along exterior bounadry $ {\bm{c}}_{0,j} $ (the radius is $ \beta_{j} $) for rotation normalization, and generally can take the coordinate $ w_{j}^{\beta} = \beta $ for simplicity; $ P_{j,k} $ and $ Q_{j,k} $ denote {\emph{charges}} of Charge Simulation Method, and are real variables to be determined by Eq. (\ref{eqa:5}); $ U_{j,k} $ and $ V_{j,k} $ denote {\emph{charge points}} of Charge Simulation Method, and are complex coefficients determined by Eq. (\ref{eqa:3}); $ p_{j,k} $ and $ q_{j,k} $ are {\emph{charges}} of Complex Dipole Simulation Method, and are complex variables to be determined by Eq. (\ref{eqa:6}); $ \eta_{j,k} $ and $ \mu_{j,k} $ are {\emph{charge points}} of Complex Dipole Simulation Method, and are complex coefficients determined by Eq. (\ref{eqa:4}); $ N_{0} $ and $ N_{j} $ denote quantities of charge points along exterior boundary $ {\bm{c}}_{0,j}^{w} $ (or $ {\bm{c}}_{0,j} $) and interior boundary $ {\bm{c}}_{j}^{w} $ (or $ {\bm{c}}_{j} $), respectively.

Via Eq. (\ref{eqa:2}), the free arc segment $ {\bm{c}}_{0f,j}^{w} $ and constrained arc segment $ {\bm{c}}_{0c,j}^{w} $ of exterior boundary $ {\bm{c}}_{0f}^{w} $ in the interval mapping plane $ w_{j} = \Re w_{j} + {\rm{i}} \Im w_{j} $ are bidirectinally mapped onto free segment $ {\bm{c}}_{0f,j} $ and constrained arc segment $ {\bm{c}}_{0c,j} $ of unit exterior boundary $ {\bm{c}}_{0,j} = {\bm{c}}_{0f,j} \cup {\bm{c}}_{0c,j} $ in the mapping plane $ \zeta_{j} = \rho_{j} \cdot {\rm{e}}^{{\rm{i}}\theta_{j}} $, respectively. The joint points $ t_{1,j}^{w} $ and $ t_{2,j}^{w} $ connecting arc segments $ {\bm{c}}_{0f,j}^{w} $ and $ {\bm{c}}_{0c,j}^{w} $ in the interval mapping plane $ w = \Re w + {\rm{i}} \Im w $ are bidirectionally mapped onto corresponding joint points $ t_{1,j} $ and $ t_{2,j} $ connecting unit arc segments $ {\bm{c}}_{0f,j} $ and $ {\bm{c}}_{0c,j} $ in the mapping plane $ \zeta_{j} = \rho_{j} \cdot {\rm{e}}^{{\rm{i}}\theta_{j}} $, respectively. The possibly noncircular interior boundary $ {\bm{c}}_{j}^{w} $ in the interval mapping plane $ w = \Re w + {\rm{i}}\Im w $ is bidirectionally mapped onto the circular interior boundary $ {\bm{c}}_{j} $ of radius $ \alpha_{j} $. The interior and exterior regions of boundary $ {\bm{c}}_{0,j} $ are denoted by $ {\bm{\omega}}_{0,j}^{+} $ and $ {\bm{\omega}}_{0,j}^{-} $, respectively. Figs. \ref{fig:a1}b and \ref{fig:a1}c graphically and schematically show the subsequent bidirectional conformal mapping in Eq. (\ref{eqa:2}) for the example in Fig. \ref{fig:a1}a.

The charge points of Charge Simulation Method can be given as
\begin{subequations}
  \label{eqa:3}
  \begin{equation}
    \label{eqa:3a}
    \left\{
      \begin{aligned}
        & U_{j,k} = w_{j,k}^{0} + K_{j}^{0} \cdot H_{j,k}^{0} \cdot \exp\left( {\rm{i}}\varTheta_{j,k}^{0} \right) \\
        & H_{j,k}^{0} = \frac{1}{2}\left( \left| w_{j,k+1}^{0} - w_{j,k}^{0} \right| + \left| w_{j,k}^{0} - w_{j,k-1}^{0} \right| \right) \\
        & \varTheta_{j,k}^{0} = \arg\left( w_{j,k+1}^{0} - w_{j,k-1}^{0} \right) - \frac{\pi}{2}
      \end{aligned}
    \right.
    , \quad k = 1,2,3,\cdots,N_{0}
  \end{equation}
  \begin{equation}
    \label{eqa:3b}
    \left\{
      \begin{aligned}
        & V_{j,k} = w_{j,k}^{1} + K_{j}^{1} \cdot H_{j,k}^{1} \cdot \exp\left( {\rm{i}}\varTheta_{j,k}^{1} \right) \\
        & H_{j,k}^{1} = \frac{1}{2}\left( \left| w_{j,k+1}^{1} - w_{j,k}^{1} \right| + \left| w_{j,k}^{1} - w_{j,k-1}^{1} \right| \right) \\
        & \varTheta_{j,k}^{1} = \arg\left( w_{j,k+1}^{1} - w_{j,k-1}^{1} \right) - \frac{\pi}{2}
      \end{aligned}
    \right.
    , \quad k = 1,2,3,\cdots,N_{j}
  \end{equation}
\end{subequations}
where $ w_{j,k}^{0} $ and $ w_{j,k}^{1} $ denote {\emph{collocation points}} of Charge Simulation Method distributed along exterior boundary $ {\bm{c}}_{0,j}^{w} $ and interior boundary $ {\bm{c}}_{j}^{w} $ of the interval mapping region $ \overline{\bm{\varOmega}}_{j}^{w} $, respectively, and the selection schemes can be seen in \ref{sec:a2}; $ K_{j}^{0} $ and $ K_{j}^{1} $ denote {\emph{assignment factors}} of Charge Simulation Method for collocation points along exterior boundary $ {\bm{c}}_{0,j}^{w} $ and interior boundary $ {\bm{c}}_{j}^{w} $, respectively; $ H_{j,k}^{0} $ and $ H_{j,k}^{1} $ denote distances between charge points and corresponding collocation points along exterior boundary $ {\bm{c}}_{0,j}^{w} $ and interior boundary $ {\bm{c}}_{j}^{w} $, respectively; $ \varTheta_{j,k}^{0} $ and $ \varTheta_{j,k}^{1} $ denote outward normal angles pointing from collocation points to corresponding charge points along exterior boundary $ {\bm{c}}_{0,j}^{w} $ and interior boundary $ {\bm{c}}_{j}^{w} $, respectively.

The charge points of Complex Dipole Simulation Method can be given as
\begin{subequations}
  \label{eqa:4}
  \begin{equation}
    \label{eqa:4a}
    \left\{
      \begin{aligned}
        & \eta_{j,k} = \zeta_{j,k}^{0} + k_{j}^{0} \cdot h_{j,k}^{0} \cdot \exp\left( {\rm{i}}\vartheta_{j,k}^{0} \right) \\
        & h_{j,k}^{0} = \frac{1}{2}\left( \left| \zeta_{j,k+1}^{0} - \zeta_{j,k}^{0} \right| + \left| \zeta_{j,k}^{0} - \zeta_{j,k-1}^{0} \right| \right) \\
        & \vartheta_{j,k}^{0} = \arg\left( \zeta_{j,k+1}^{0} - \zeta_{j,k-1}^{0} \right) - \frac{\pi}{2}
      \end{aligned}
    \right.
    , \quad k = 1,2,3,\cdots,N_{0}
  \end{equation}
  \begin{equation}
    \label{eqa:4b}
    \left\{
      \begin{aligned}
        & \mu_{j,k} = \zeta_{j,k}^{1} + k_{j}^{1} \cdot h_{j,k}^{1} \cdot \exp\left( {\rm{i}}\vartheta_{j,k}^{1} \right) \\
        & h_{j,k}^{1} = \frac{1}{2}\left( \left| \zeta_{j,k+1}^{1} - \zeta_{j,k}^{1} \right| + \left| \zeta_{j,k}^{1} - \zeta_{j,k-1}^{1} \right| \right) \\
        & \vartheta_{j,k}^{1} = \arg\left( \zeta_{j,k+1}^{1} - \zeta_{j,k-1}^{1} \right) - \frac{\pi}{2}
      \end{aligned}
    \right.
    , \quad k = 1,2,3,\cdots,N_{j}
  \end{equation}
\end{subequations}
where $ \zeta_{j,k}^{0} $ and $ \zeta_{j,k}^{1} $ denote {\emph{collocation points}} of Complex Dipole Simulation Method distributed along exterior boundary $ {\bm{c}}_{0,j} $ and interior boundary $ {\bm{c}}_{j} $ of the interval mapping region $ \overline{\bm{\omega}}_{j} $, respectively, and are determined by Eq. (\ref{eqa:7}); $ k_{j}^{0} $ and $ k_{j}^{1} $ denote {\emph{assignment factors}} of Complex Dipole Simulation Method for collocation points along exterior boundary $ {\bm{c}}_{0,j} $ and interior boundary $ {\bm{c}}_{j} $, respectively; $ h_{j,k}^{0} $ and $ h_{j,k}^{1} $ denote distances between charge points and corresponding collocation points along exterior boundary $ {\bm{c}}_{0,j} $ and interior boundary $ {\bm{c}}_{j} $, respectively; $ \vartheta_{j,k}^{0} $ and $ \vartheta_{j,k}^{1} $ denote outward normal angles pointing from collocation points to corresponding charge points along exterior boundary $ {\bm{c}}_{0,j} $ and interior boundary $ {\bm{c}}_{j} $, respectively.

With the charge points in Eq. (\ref{eqa:3}), the charges in Eq. (\ref{eqa:2a}) can be determined as
\begin{subequations}
  \label{eqa:5}
  \begin{equation}
    \label{eqa:5a}
    \sum\limits_{k=1}^{N_{0}} \ln\left| \frac{w_{j,i}^{0} - U_{j,k}}{w_{j}^{\beta} - U_{j,k}} \right| \cdot P_{j,k} + \sum\limits_{k=1}^{N_{j}} \ln\left| \frac{w_{j,i}^{0} - V_{j,k}}{w_{j}^{\beta} - V_{j,k}} \right| \cdot Q_{j,k} - \ln r_{e} = - \ln \left| \frac{w_{j,k}^{0} - w_{j}^{c}}{w_{j}^{\beta} - w_{j}^{c}} \right|, \quad i = 1,2,3,\cdots,N_{0}
  \end{equation}
  \begin{equation}
    \label{eqa:5b}
    \sum\limits_{k=1}^{N_{0}} \ln\left| \frac{w_{j,i}^{1} - U_{j,k}}{w_{j}^{\beta} - U_{j,k}} \right| \cdot P_{j,k} + \sum\limits_{k=1}^{N_{j}} \ln\left| \frac{w_{j,i}^{1} - V_{j,k}}{w_{j}^{\beta} - V_{j,k}} \right| \cdot Q_{j,k} - \ln \alpha_{j} = - \ln \left| \frac{w_{j,k}^{1} - w_{j}^{c}}{w_{j}^{\beta} - w_{j}^{c}} \right|, \quad i = 1,2,3,\cdots,N_{j}
  \end{equation}
  \begin{equation}
    \label{eqa:5c}
    \sum\limits_{k=1}^{N_{0}} P_{j,k} = -1
  \end{equation}
  \begin{equation}
    \label{eqa:5d}
    \sum\limits_{k=1}^{N_{j}} Q_{j,k} = 0
  \end{equation}
\end{subequations}
Eq. (\ref{eqa:5}) contains $ N_{0}+N_{j}+2 $ to-be-determined real variables and real linear equations to form simultaneous real linear system, and unique solution can be obtained.

With the charge points in Eq. (\ref{eqa:4}), the charges in Eq. (\ref{eqa:2b}) can be determined as
\begin{subequations}
  \label{eqa:6}
  \begin{equation}
    \label{eqa:6a}
    \sum\limits_{k=1}^{N_{0}} \frac{p_{j,k}}{\zeta_{j,i}^{0}-\eta_{j,k}} + \sum\limits_{k=1}^{N_{j}} \frac{q_{j,k}}{\zeta_{j,i}^{0}-\mu_{j,k}} = w_{j,i}^{0}, \quad i = 1,2,3,\cdots,N_{0}
  \end{equation}
  \begin{equation}
    \label{eqa:6b}
    \sum\limits_{k=1}^{N_{0}} \frac{p_{j,k}}{\zeta_{j,i}^{1}-\eta_{j,k}} + \sum\limits_{k=1}^{N_{j}} \frac{q_{j,k}}{\zeta_{j,i}^{1}-\mu_{j,k}} = w_{j,i}^{1}, \quad i = 1,2,3,\cdots,N_{j}
  \end{equation}
\end{subequations}
where
\begin{subequations}
  \label{eqa:7}
  \begin{equation}
    \label{eqa:7a}
    \zeta_{j,i}^{0} = \zeta_{j}(w_{j,i}^{0}), \quad i = 1,2,3,\cdots,N_{0}
  \end{equation}
  \begin{equation}
    \label{eqa:7b}
    \zeta_{j,i}^{1} = \zeta_{j}(w_{j,i}^{1}), \quad i = 1,2,3,\cdots,N_{j}
  \end{equation}
\end{subequations}
Eq. (\ref{eqa:7}) is computed using Eq. (\ref{eqa:2a'}). Eq. (\ref{eqa:6}) contains $ N_{0}+N_{j} $ to-be-determined complex variables and complex linear equations to form simultaneous complex linear system, and unique solution can be obtained.

In summary, the bidirectional conformal mapping in Eq. (\ref{eqa:2}) can be determined by the linear systems in Eqs. (\ref{eqa:5}) and (\ref{eqa:6}) with auxiliary supports of Eqs. (\ref{eqa:3}), (\ref{eqa:4}), and (\ref{eqa:7}), except for the before-hand determination of collocation points $ w_{j,k}^{0} $ and $ w_{j,k}^{1} $, which are given in \ref{sec:a2} with certain numerical principles.

\subsection{Numerical schemes}
\label{sec:a2}

\begin{figure}[htb]
  \centering
  \includegraphics[width = \textwidth]{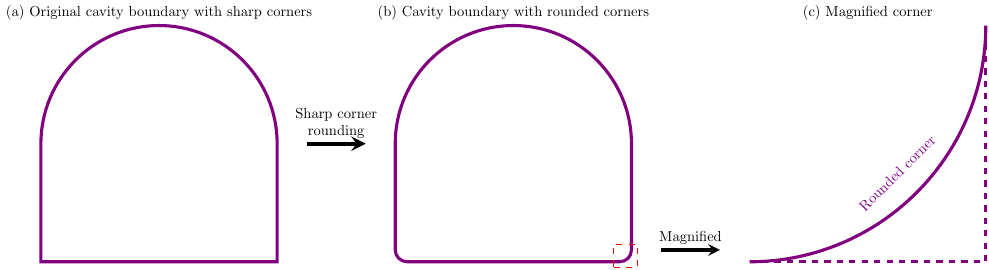}
  \caption{Rounding of sharp corners -- schematic diagram for the cavity boundary in Fig. \ref{fig:a1}a}
  \label{fig:a2}
\end{figure}

The collocation point selection of $ w_{j,k}^{0} $ and $ w_{j,k}^{1} $ in Eq. (\ref{eqa:3}) should be suitable to guarantee the analyticity and accuracy of the backward conformal mapping in Eq. (\ref{eqa:2b}). The collocation points along boundaries $ {\bm{c}}_{0,j}^{w} $ and $ {\bm{c}}_{j}^{w} $ of the interval mapping annulus $ \overline{\bm{\varOmega}}_{j}^{w} $ should be distributed as uniformly as possible, and sharp corners should be eliminated to avoid boundary singularities \cite{hough1981use,hough1983integral}. Therefore, we should first ensure that no sharp corners exist along the cavity boundary $ {\bm{C}}_{j} $ in the physical geomaterial region $ \overline{\bm{\varOmega}}_{j} $. In Figs. \ref{fig:2} and \ref{fig:3} and \ref{fig:a1}a, the sharp corners of cavity boundaries for sequential shallow tunnelling are kept merely due to convenience of schematic plotting. In practical tunnel engineering, sharp corners are generally rounded on purpose to avoid stress concentration. Therefore, the sharp corners in this paper do not exist, and Fig. \ref{fig:a2} graphically shows the replacement of rounded corners with sharp corners for the example of Fig. \ref{fig:a1}a.

Therefore, the rounded corners guarantee the continuity and analyticity of the cavity boundary $ {\bm{C}}_{j} $ in the physical plane $ z = x+{\rm{i}}y $. Subsequently, via the first-step conformal mapping in Eq. (\ref{eqa:1}), the exterior circular boundary $ {\bm{c}}_{0,j}^{w} $ (mapped from the straight ground surface $ {\bm{C}}_{0} $) and interior boundary $ {\bm{c}}_{j}^{w} $ in the inteval mapping annulus $ \overline{\bm{\varOmega}}_{j}^{w} $ (possibly not circular) would contain no sharp corners, so that the continuity of the outward normal angles $ \varTheta_{j,k}^{0} $ and $ \varTheta_{j,k}^{1} $ in Eq. (\ref{eqa:3}) would be potentially guaranteed by mandatorily setting the following numerical principle as
\begin{subequations}
  \label{eqa:8}
  \begin{equation}
    \label{eqa:8a}
    \max\left[ \max\limits_{1 \leq k \leq N_{0}}\left| \arg\left( w_{j,k+1}^{0} - w_{j,k}^{0} \right) \right|, \max\limits_{1 \leq k \leq N_{j}}\left| \arg\left( w_{j,k+1}^{1} - w_{j,k}^{1} \right) \right|  \right] \leq 10^{\circ}
  \end{equation}
  Moreover, the quantity of collocation points should be large enough, thus, the following numerical principle is mandatorily set as
  \begin{equation}
    \label{eqa:8b}
    \max\left[ \max\limits_{1 \leq k \leq N_{0}} \frac{| w_{j,k+1}^{0} - w_{j,k}^{0} |}{\sum\limits_{l=1}^{N_{0}}| w_{j,l}^{0} - w_{j,l}^{0} |}, \max\limits_{1 \leq k \leq N_{j}} \frac{| w_{j,k+1}^{1} - w_{j,k}^{1} |}{\sum\limits_{l=1}^{N_{j}}| w_{j,l}^{1} - w_{j,l}^{1} |} \right] \leq 10^{-2}
  \end{equation}
\end{subequations}
Eq. (\ref{eqa:8}) is based on the numerical principles proposed in our previous study \cite{lin2024charge-simulation}, and controls the varying angle and distance between adjacent collocation points along boundaries $ {\bm{c}}_{0,j}^{w} $ and $ {\bm{c}}_{j}^{w} $.

Eq. (\ref{eqa:2a}) should be slightly modified after solving the charges by Eq. (\ref{eqa:5}) to guarantee the single-valuedness of arguments without changing the principal values of the computation results \cite{okano2003numerical,lin2024charge-simulation} as
\begin{equation}
  \label{eqa:2a'}
  \tag{A.2a'}
  \zeta_{j} = \zeta_{j}(w_{j}) = \frac{w_{j}-w_{j}^{c}}{w_{j}^{\beta}-w_{j}^{c}} \cdot \exp\left[
    \begin{aligned}
      &\sum\limits_{k=1}^{N_{0}} P_{j,k} \left( \ln\frac{w_{j}-U_{j,k}}{w_{j,0}-U_{j,k}} - \ln\frac{w_{j}^{\beta}-U_{j,k}}{w_{j,0}-U_{j,k}} \right) \\
      + & \sum\limits_{k=1}^{N_{j}} Q_{j,k} \left( \ln\frac{w_{j}-V_{j,k}}{w_{j}-w_{j}^{c}} - \ln\frac{w_{j}^{\beta}-V_{j,k}}{w_{j}^{\beta}-w_{j}^{c}} \right)
    \end{aligned}
  \right]
\end{equation}
where $ w_{j,0} $ denote an arbitrary point within interval mapping region $ \overline{\bm{\varOmega}}_{j}^{w} $.

%%===========================================================================================%%
%% If you are submitting to one of the Nature Portfolio journals, using the eJP submission   %%
%% system, please include the references within the manuscript file itself. You may do this  %%
%% by copying the reference list from your .bbl file, paste it into the main manuscript .tex %%
%% file, and delete the associated \verb+\bibliography+ commands.                            %%
%%===========================================================================================%%

\bibliography{biblio}% common bib file
%% if required, the content of .bbl file can be included here once bbl is generated
%%\input sn-article.bbl

\end{document}